\documentclass[a4paper, reqno, 10pt]{amsart}
\usepackage[foot]{amsaddr}

\usepackage[T1]{fontenc}
\usepackage[utf8]{inputenc}
\usepackage[english]{babel}
\usepackage{amsmath, amssymb, amsthm}
\usepackage[scr=rsfso]{mathalpha}
\usepackage{enumitem} 
\usepackage{multicol} 
\usepackage[framemethod=TikZ]{mdframed}
\usepackage{tikz-cd}
\usepackage{caption} 
\usepackage{subcaption} 
\usepackage{cite}
\usepackage{float}
\usepackage{stackrel}
\usepackage{mathtools} 
\usepackage{extarrows} 
\usepackage{bbm} 
\usepackage[export]{adjustbox} 
\usepackage{etoolbox} 
\usepackage{xurl}

\usepackage[linktoc=all]{hyperref} 
\usepackage[capitalise]{cleveref} 
\usepackage{xcolor} 

\hypersetup{
	colorlinks,
	linkcolor={red!85!black}, 
	citecolor={green!55!black},
	urlcolor={blue!85!black}
}

\usepackage{algpseudocode}
\usepackage{algorithm}

\usepackage{subfiles} 

\makeatletter
\def\paragraph{\@startsection{paragraph}{4}%
	\z@\z@{-\fontdimen2\font}%
	{\normalfont\bfseries}}
\makeatother

\title{Relative field theories via relative dualizability}
\author[C.~Scheimbauer]{Claudia Scheimbauer}
\address{TU M{\"{u}}nchen, Boltzmannstra\ss e 3, 85748 Garching, Germany, 
ORCiD:0000-0002-6927-8348
}
\email{scheimbauer@tum.de (corresponding author)}

\author[T.~Stempfhuber]{Thomas Stempfhuber}
\address{Universit{\"{a}}t Hamburg, Bundesstra\ss e 55,
	20146 Hamburg, Germany}
\email{stempfhuber.thomas@gmail.com}
\date{\today}

\makeatletter
\newtheorem*{rep@theorem}{\rep@title}
\newcommand{\newreptheorem}[2]{%
\newenvironment{rep#1}[1]{%
 \def\rep@title{#2 \ref{##1}}%
 \begin{rep@theorem}}%
 {\end{rep@theorem}}}
\makeatother

\newtheorem{theo}{Theorem}[section]
\newreptheorem{theo}{Theorem}
\newtheorem{coro}[theo]{Corollary}
\newtheorem{lema}[theo]{Lemma}
\newtheorem{prop}[theo]{Proposition}

\theoremstyle{definition}
\newtheorem{defi}[theo]{Definition}

\newcommand{\exampleend}{\hfill \ensuremath{\diamondsuit}}
\theoremstyle{remark}
\newtheorem{rema}[theo]{Remark}
\newtheorem{exam}[theo]{Example}
\newtheorem*{nota}{Notation}

\DeclareTextFontCommand{\emph}{\bfseries} 

\DeclareMathOperator{\Hom}{Hom} 
\DeclareMathOperator{\id}{id} 
\DeclareMathOperator{\ho}{h} 
\DeclareMathOperator{\fd}{fd} 

\newcommand{\op}{\mathrm{op}}

\newcommand{\Cat}{\textbf{Cat}}
\newcommand{\Vect}{\textbf{Vect}}

\newcommand{\Bord}{\textbf{Bord}}
\newcommand{\Alg}[1]{\textbf{Alg}_{#1}}
\newcommand{\LinCat}{\textbf{LinCat}}

\newcommand{\TFT}[1]{\mathscr{#1}}

\newcommand\opcat[1]{{#1}^{\mathrm{op}}}
\newcommand\cocat[1]{{#1}^{\mathrm{co}}}
\newcommand\coopcat[1]{{#1}^{\mathrm{coop}}}

\newcommand{\Cri}{\ca C^{\rightarrow}}
\newcommand{\Cdo}{\ca C^{\downarrow}}
\newcommand{\CaRight}[1]{{#1}^{\rightarrow}}

\DeclareMathOperator{\even}{even}
\DeclareMathOperator{\odd}{odd}

\newcommand{\tree}[1]{{\mathsf {T_{#1}} }}
\newcommand{\treeSet}[1]{\mathcal{T}_{#1} }
\newcommand{\treeEqSet}[1]{\mathscr{T}_{#1} }

\newcommand{\ca}[1]{{\mathcal{#1} }} 
\newcommand{\infN}[1]{(\infty , #1)\text{-}} 
\newcommand{\te}{\otimes} 
\newcommand{\uobj}{\mathbbm{1}} 
\newcommand{\adj}{\dashv} 

\newcommand{\RLarrows}{\mathrel{\substack{\longrightarrow \\[-.6ex] \longleftarrow}}}

\newcommand{\LongDoubleRLarrows}{\substack{\xRightarrow{\qquad} \\[-.9ex] \xLeftarrow{\qquad} }}

\makeatletter
\newcommand{\xRrightarrow}[2][]{\ext@arrow 0359\Rrightarrowfill@{#1}{#2}}
\newcommand{\Rrightarrowfill@}{\arrowfill@\equiv\equiv\Rrightarrow}
\makeatother


\begin{document}


\begin{abstract} \noindent
We investigate relative versions of dualizability designed for relative versions of topological field theories (TFTs), also called twisted TFTs, or quiche TFTs in the context of symmetries.
In even dimensions we show an equivalence between lax and oplax fully extended framed relative topological field theories
valued in an $\infN N$category in terms of adjunctibility.
Motivated by this, we systematically investigate higher adjunctibility conditions and their implications for relative TFTs.
Summarizing we arrive at the conclusion that oplax relative TFTs is the notion of choice.
Finally, for fun we explore a tree version of adjunctibility and compute the number of equivalence classes thereof.
\end{abstract}

\maketitle


{
	\hypersetup{hidelinks} 
	\tableofcontents
}

\thispagestyle{empty}

\setcounter{page}{1}



\section{Introduction} \label{sec:intro}

The notion of a functor having an adjoint is a property describing many phenomena in all areas of mathematics. Replacing the 2-category of categories, functors, and natural transformations by an arbitrary bicategory, or, even more generally, an $(\infty,n)$-category, we obtain a notion which generalizes --  simultaneously -- finiteness conditions for a module (finitely generated and projective), $k$-handles of a manifold and handle cancellation, duals of vector spaces and perfect complexes, and smooth and proper dg algebras or dg categories.

In fact, we have a hierarchy of properties, which we will investigate in this article. We start with a 1-morphism in our higher category $\mathcal{C}$ and ask for this 1-morphism to have a left adjoint (or to have a right adjoint). This adjoint comes together with a unit and counit 2-morphism exhibiting adjunctibility. Now we can ask for these 2-morphisms themselves to have a left adjoint (or to have a right adjoint). In an $(\infty,N)$-category, we can repeat this procedure $n$ times, arriving at {\em $n$-times left (or right) adjunctibility}.

In this article, our first goal was to demonstrate that for even $n$, opting for always either left or right adjoints at every level leads to the same notion. This result is established through an interchange lemma generalizing such an interchange property stated in Lurie's seminal article on the Cobordism Hypothesis~\cite{lurie2009classification}.

While this result is a purely higher categorical result, our motivation for proving this  comes from the classification of a {\em relative} version of fully extended topological field theories \cite{FT2014relative, stolz2011supersymmetric}.
Such a relative notion of field theory was  first introduced by Stolz and Teichner in \cite{stolz2011supersymmetric} as \emph{twisted field theories} to capture the behavior of Segal's conformal field theories \cite{CFT} and anomalies: rather than assigning a number to a closed top dimensional manifold viewed as an endomorphism of the ground field $k$, a twisted field theory should capture choosing elements in a vector space or line, viewed as a morphism $k\to V$. Here the vector space may vary depending on the choice of manifold, hence is a once-categorified field theory itself.
Freed and Teleman's \emph{relative field theories} \cite{FT2014relative} are meant to capture the same idea, namely, they should be a ``homomorphism'' of field theories $\uobj\to \alpha$, where $\alpha$ often (but not always) is the truncation of a field theory of one dimension higher.

The mathematical framework suggested in \cite{stolz2011supersymmetric} to capture this behavior is that of a \emph{symmetric monoidal lax or oplax natural transformation}
\begin{equation} \label{eq:TwistedTQFT}
\begin{tikzcd}
   &                                                                     & {} \arrow[dd, "\Longrightarrow" rotate=-90, phantom, shift right] &           \\
& \Bord_n \arrow[rr, "\uobj", bend left] \arrow[rr, "T"', bend right] &                                       & \ \ca C \ \\
   &                                                                     & {}                                    &          
\end{tikzcd}
\end{equation}
between symmetric monoidal functors. Here $\ca C$ is a bicategory, $\uobj$ is the trivial functor, and $T$ is
called twist.
In \cite{JF_CS} this notion was extended to the higher categorical setting. Moreover, using the Cobordism Hypothesis, a classification of twisted fully extended topological field theories was given precisely in terms of left (or right)  $n$-times adjunctibility.
Our main result shows that in even dimensions, these flavors occur simultaneously; while in odd dimensions, they differ.

Finally, our interchange result can be generalized: so far we have chosen either left at every stage, or right at every stage. Instead, we may choose to alternate in a chosen pattern.
Our interchange lemma shows that in any dimension there are exactly two equivalence classes of such notions of adjunctibility. In even dimensions left and  right $n$-times-adjunctibility happen to capture just one of the classes, whereas in odd dimensions, they lie in different classes.
One may ask for field theoretic interpretations of all of these notions, and we discuss this in detail.
In particular, we answer the question which notion of adjunctibility is the ``correct'' one for relative field theories: we should ask for an \emph{oplax} natural transformation\footnote{We explain two desiderata for relative field theories below in \eqref{eqn:element_in_vsp_intro} and \eqref{eqn:deloopTFTintro}. We prove that they cannot be simultaneously satisfied, and priority is given to \eqref{eqn:element_in_vsp_intro}. Moreover, we prove that the other choices can be related to the oplax case by taking various opposite categories of the target.}.

Before delving into the details of our categorical results, let us mention that these relative versions of (T)FTs have recently been ubiquitous: they appear as ``quiche'' field theories\footnote{Freed--Moore--Teleman in \cite{FMT} introduced the term {\em quiche} as being half of a {\em sandwich}. We would like to point out here that a quiche is not half of a sandwich as was already observed by the authors themselves. In German the word {\em Butterbrot} \url{https://en.wikipedia.org/wiki/Butterbrot} is precisely half of a sandwich, namely a piece of bread with a topping of butter (and in some regions some extra stuff). We refrain from introducing yet another different name.\label{foot_Butterbrot}}
 in the context of (possibly non-invertible) symmetries in \cite{FMT}, see in particular Remarks 2.7 (3) and 3.3 (8). Using the framework developed in \cite{JF_CS} which is central here,  various such examples have recently been worked out, e.g.~in \cite{haioun2023unit, vandyke2023projective, vandyke2023symmetries, kinnear2024nonssCY}.

Finally, let us mention that in \cite{StewartThesis} the connection of this notion of relative field theory using a categorical framework to a geometric definition of relative field theory using a variation of a bordism category with defects or boundaries was proven. In particular, in the framed case different choices of how the framing should restrict along the defect or boundary correspond precisely to the different choices of adjunctibility presented in this article.

\section*{Overview of our higher categorical results}

We now give a more detailed overview of our higher categorical results.

Our first main result, which was the motivation for this article,
is that in an $\infN {N}$category $\ca C$ and for even $n>0$, $n$-times left adjunctibility and $n$-times right adjunctibility of a $k$-morphism coincide. This follows from an interchange result,  \cref{Lemma:MainLemma}.
\begin{reptheo}{Coro:evenDimensions}
    Let $n\geq 2$ be even. Then a $k$-morphism $f$ is $n$-times left adjunctible
    if and only if $f$ is $n$-times right adjunctible.
\end{reptheo}

In the odd case, this is not true. In fact, if a $k$-morphism is both $n$-times left and $n$-times right adjunctible, it already is $n$-adjunctible (\cref{Theo:fromMixedToNadjunctibility} ).

One may now ask whether in even dimensions there is a different notion of adjunctibility which is not equivalent to either purely right or purely left adjunctibility.
Indeed, rather than choosing the existence of right adjoints of (co)units at each level (or always left), we can use a general ``mixed'' sequence of ``right'' and ``left'' to define a notion of  ``mixed'' adjunctibility. We call such a sequence a \emph{dexterity function} $a^n$ and adjunctibility based thereupon $a^n$-adjunctibility. Given two such functions $a^n, b^n: \{ 1, 2, \ldots, n \} \rightrightarrows \{ L, R \}$, their common adjunctibility depends on the parity of 
\begin{displaymath}
	p(a^n, b^n) \ \coloneqq \ \	| (a^n)^{-1}(L) | \, - \, | (b^n)^{-1}(L) | \,.
	\end{displaymath}
The Theorem above then generalizes to
\begin{theo} [\cref{Prop:AnReducesToTwoClasses}]
If $p(a^n, b^n)$ is even,
    then  $f$ is $a^n$-adjunctible if and only if $f$ is $b^n$-adjunctible.
\end{theo}
In fact, the $2^n$ a priori different notions of mixed adjunctibility reduce to just two: for every $n$ there are \emph{exactly} two equivalence classes of ``mixed'' adjunctibility, where we call two dexterity functions $a^n, b^n$ and their notions of mixed adjunctibility  equivalent if the following holds: a morphism is $a^n$-adjunctible if and only if it is $b^n$-adjunctible. We prove that this is the case if and only their parity $p(a^n, b^n)$ is even. This is \cref{cor:parity_implies_even/odd} and \cref{rem:exactly_two_classes}.

If their parity is odd, we obtain a stronger condition.
\begin{theo} [\cref{Prop:AnReducesToTwoClasses} and \cref{Theo:fromMixedToNadjunctibility} ]
If  $p(a^n, b^n)$ is odd and   $f$  is both $a^n$-adjunctible and $b^n$-adjunctible, then $f$ is $n$-adjunctible.
\end{theo}
Along the way we prove some cute results useful to proving various adjunctibility results.

For a bicategory, if we reverse only 1-morphisms, left adjoints become right adjoints and vice versa. In \cref{sec:opposite} we consider various opposite categories, reversing only some of the morphisms, in relation to mixed adjunctibility. More precisely, for two dexterity functions $a^n, b^n$ we can precisely say for which $j$ we have to reverse $j$-morphisms in $\ca C$ so that $a^n$-adjunctibility translates to $b^n$-adjunctibility.

Our main result  suggests that for $n$ odd, there should be
a morphism $f$ which is $n$-times left but not right adjunctible.
In \cref{sec:sec2}, \cref{ex:3-times-left-not-right} we construct an example illustrating this for $n=3$ in a higher Morita category.

\section*{Application to relative and twisted topological field theories}

Relative versions of field theories should satisfy two desiderata, which we elaborate on in \cref{sec:applications}.  The first, already mentioned above, is that it should for top-dimensional closed manifolds, we should choose an element in a vector space (for instance a line), \begin{equation}\label{eqn:element_in_vsp_intro}
\uobj \to T(M) \,.
\end{equation}
The second is that if $T\equiv \uobj$, we should recover the usual notion of ``absolute'' field theory valued in the looping $\Omega\mathcal{C}$,
\begin{equation}\label{eqn:deloopTFTintro}
\TFT{Z}\colon \Bord_n^{fr} \to \Omega\mathcal{C}\,.
\end{equation}
The framework suggested in \cite{stolz2011supersymmetric} using symmetric monoidal lax or oplax natural transformations unfortunately does not satisfy both at the same time: As was shown in \cite{JF_CS},  oplax field theories satisfy the former, but not the latter; whereas lax field theories satisfy the latter, but not the former.

The main goal of our categorical results was to compare the oplax and lax variant of relative field theories, and to see whether we can modify the framework to obtain both  \eqref{eqn:element_in_vsp_intro} and \eqref{eqn:deloopTFTintro}. Moreover, Freed--Teleman in \cite{FT2014relative} did not make precise which framework they would like to use, and for instance in \cite{FreedTelemanGapped} suggest a stronger condition in terms of adjunctibility. In this article we study the various options for the framework and compare them to this latter stronger condition.

Recall that an oplax (respectively lax) natural transformation is a functor
\[\ca B \longrightarrow \Cri  \quad \mbox{respectively} \quad \ca B \longrightarrow \Cdo \,,\]
where $\Cri$ and $\Cdo$ are certain arrow higher categories governing the desired diagrammatics.
Hence an oplax (respectively lax) relative field theory is a symmetric monoidal functor
\[\Bord_n \longrightarrow \Cri  \quad \mbox{respectively} \quad \Bord_n \longrightarrow \Cdo \,.\]

In the fully extended framed case (using $\Bord_n^{fr}$) and under the assumption of the \emph{Cobordism Hypothesis} from \cite{lurie2009classification}, lax and oplax 
twisted topological field theories are classified by $n$-dualizable objects in the
target categories $\Cdo$ and $\Cri$.
These were in turn computed by Johnson-Freyd and the first author.
\begin{reptheo}{theo:laxOplaxDualizable}[Theorem 7.6 in \cite{JF_CS}]
  An object $f:X \to Y$ in $\Cdo$ is $n$-dualizable
  if and only if $X, Y$ are $n$-dualizable and
  $f$ as a morphism in $\ca C$ is $n$-times left adjunctible.
  An object $f:X \to Y$ in $\Cri$ is $n$-dualizable
  if and only if $X, Y$ are $n$-dualizable and
  $f$ as a morphism in $\ca C$ is $n$-times right adjunctible.
\end{reptheo}

Together with the Cobordism Hypothesis, this Theorem reduces the comparison of lax and oplax fully extended framed
twisted topological field theories to the categorical comparison of $n$-times left and
right adjunctibility, which is our main categorical result \cref{Coro:evenDimensions}.

This leads to the question whether, for any given dexterity function $a^n$,  we can define arrow categories ${\ca C}^{a^n}$ whose $n$-dualizable objects are characterized by $a^n$-adjunctible morphisms as in \cref{theo:laxOplaxDualizable}. More importantly, can we choose a dexterity function $a^n$ so that we resolve our desiderata \eqref{eqn:element_in_vsp_intro} and \eqref{eqn:deloopTFTintro}?

We answer this question negatively. To see this, we first prove that for desideratum \eqref{eqn:element_in_vsp_intro}, choosing \emph{oplax} is the only possible choice. This follows from 
\begin{reptheo}{cor:a-lax_dualizability}
	Let $\ca C$ be a symmetric monoidal $\infN N$category, $f:X \to Y$ a $1$-morphism in $\ca C$, and $a^n$ a dexterity function. The following  are equivalent:
	\begin{enumerate}
	\item When viewed as an object in ${\ca C}^{a^n}$, $f$ is $n$-dualizable.
	\item The objects $X$ and $Y$ are $n$-dualizable and $f$ is $n$-times right adjunctible in~${\ca C}^{\op_{a^n}}$.
	\item The objects $X$ and $Y$ are $n$-dualizable and $f$ is $a^n$-adjunctible in~$\ca C$.
	\end{enumerate}
\end{reptheo}
Here ${\ca C}^{\op_{a^n}}$ is a a certain opposite category of $\ca C$, reversing $j$-morphisms for certain~$j$ depending on $a^n$.

Our conclusion now is the following: we \emph{could} define an $a^n$-lax twisted field theory to be a symmetric monoidal ``$a^n$-lax natural transformation'', i.e.~a symmetric monoidal functor
\[\Bord_n \longrightarrow {\ca C}^{a^n},\]
and \cref{cor:a-lax_dualizability} characterizes the framed fully extended framed ones in terms of $a^n$-adjunctibility.
However, only choosing \emph{lax} satisfies desideratum \eqref{eqn:deloopTFTintro}, so 
this does not resolve the dichotomy between \eqref{eqn:element_in_vsp} and \eqref{eqn:deloopTFT}:

Finally, we investigate the above mentioned stronger adjunctibility condition appearing in \cite{FreedTelemanGapped}. Teleman suggested to us that this condition might appear when looking at dualizable field theories. We confirm this in \cref{prop:dualizable_relative}, but the converse direction  is not true. Instead, we give a characterization in terms of a stronger condition in \cref{coro:characterization_dualizable_twisted}.

\section*{A fun variant: trees of adjunctibility}

In the last section \cref{sec:moreGeneral}, we play with the higher adjunctibility conditions. So far, for fixed $j$, we required for each (co)unit $j$-morphism the existence of the same sided adjoint. We generalize this to asking for each unit and counit individually whether it should have a left or right adjoint. Which sided adjoint we require is now recorded in a \emph{dexterity tree} (rather than a dexterity function). For instance, for a morphism $f$ the tree
\begin{displaymath}
		\vcenter{\hbox{\begin{tikzpicture}
					\tikzstyle{level 1}=[level distance=6mm,sibling distance=12mm]
					\tikzstyle{level 2}=[level distance=6mm,sibling distance=6mm]			
					\node {$R$}
					child {node {$L$}
						child {node {$R$}}
						child {node {$R$}}
					}
					child {node {$R$}
						child {node {$L$}}
						child {node {$L$}}
					};
		\end{tikzpicture}}}
\end{displaymath}
records the existence of
\begin{displaymath}
	\vcenter{\hbox{\begin{tikzpicture}
				\tikzstyle{level 1}=[level distance=8mm,sibling distance=60mm]
				\tikzstyle{level 2}=[level distance=8mm,sibling distance=30mm]	
				\node {$(f \adj f^R, u, c)$}
				child {node {$(u^L \adj u, u_u, c_u)$}
						child {node {$(u_u \adj u_u^R, u_{u_u}, c_{u_u})$}}
						child {node {$(c_u \adj u_u^R, u_{c_u}, c_{c_u})$}}
				}
				child {node {$(c \adj c^R, u_c, c_c)$}
						child {node {$(u_c^L \adj u_c, u_{u_c}, c_{u_c})$}}
						child {node {$(c_c^L \adj c_c, u_{c_c}, c_{c_c})$ \,.}}
				};
				
	\end{tikzpicture}}}
\end{displaymath}

The interchange result \cref{Lemma:MainLemma} results in a equivalence relation $t^n\sim s^n$ of such trees, which amounts to a morphism being $t^n$-adjunctible if and only if it is $s^n$-adjunctible, see \cref{Theo:TreeAdjunctibilityEquivalence}.

Let $\treeEqSet{n} := \treeSet{n} / \sim$ be the set of equivalence classes of dexterity trees.
\begin{reptheo}{Theo:TreeTheorem}
We have
	\begin{displaymath}
		|\treeEqSet{1}| = 2 
		\quad \text{and} \quad
		|\treeEqSet{n}| = |\treeEqSet{n-1}|^2 + 2^{2^{n-1} -1 }
		\quad \text{for} \quad n \geq 2.
	\end{displaymath}
which can be found in the OEIS as
	\href{http://oeis.org/A332757}{A332757} \cite{oeis}. 
	It also describes  the number of involutions in the n-fold iterated wreath product of $\mathbb{Z} / 2 \mathbb{Z}$ or the number of involutory automorphisms of $\tree {n+1}$.
	The first five terms are 
	$2, 6, 44, 2064$ and  $4292864$.
\end{reptheo}

\section*{Overview of the article}
In \cref{sec:introDual} we recall the notions of duals and adjoints in higher categories and $n$-dualizability.

In \cref{sec:VariantsDual} we first recall the notions of higher left and right adjunctibility and reprove the Interchange Lemma (\cref{Lemma:MainLemma}).
Then we generalize these notions to ``mixed'' versions of higher adjunctibility and prove our main theorem (\cref{Prop:AnReducesToTwoClasses}), which together with \cref{rem:exactly_two_classes} shows that we have precisely two independent notions of mixed adjunctibility. In even dimensions, higher left and right adjunctibility reduce to the same case (\cref{Coro:evenDimensions}). We continue with useful tools for proving adjunctibility and discuss various opposite categories and mixed adjunctibility therein.

We construct various examples of mixed adjunctible morphisms in higher Morita categories in \cref{sec:sec2}.

In  \cref{sec:applications} we discuss  applications of our results to relative field theories, which was the motivation for this article. 

Finally, in  \cref{sec:moreGeneral} we consider even more general notions of higher adjunctibility, depending on binary trees and compute the number of different notions that appear.

In the Appendix we include some basic lemmas about adjunctions in \cref{sec:appAdjoints}, recall and reprove a reduction of $n$-dualizability conditions in  \cref{sec:appDual}, and include full proofs of dualizability statements in the Morita bicategory in \cref{sec:appMorita}.

\section*{Acknowledgements}
We are grateful to Dan Freed, Stephan Stolz, Peter Teichner, and Constantin Teleman for prompting us to carefully analyze the adjunctibility conditions appearing in the context of relative field theories, and for many discussions on related topics.
We thank William Stewart for catching a mistake in the appendix.
CS was supported by the SFB 1085: Higher Invariants from the Deutsche Forschungsgemeinschaft
(DFG) and the Simons Collaboration on Global Categorical Symmetries (1013836).
TS was funded by the Deutsche Forschungsgemeinschaft (DFG) through the SFB 1624: 
Higher structures, moduli spaces and integrability - 506632645.

\section{Higher duals and adjoints} \label{sec:introDual}

In this section we review the terminology of adjoints and duals as well as
their higher categorical analogues. For more details, see \cite{lurie2009classification} and \cite{JF_CS}.

\subsection{Adjoints in bicategories}

\begin{defi} \label{Def:Adjunction}
 Let $\ca C$ be a bicategory and $A, B \in \ca C$. An \emph{adjunction} between two 1-morphisms 
 \begin{displaymath}
   l: A \to B, \qquad r: B \to A,
 \end{displaymath}
 consists of a pair of $2$-morphisms
 \begin{displaymath}
   u: \id_A \Rightarrow r \circ l, \qquad c: l \circ r \Rightarrow \id_B,
 \end{displaymath}
 satisfying the \emph{zig-zag identities}
 \begin{align} 
   & ( \ l \ \xRightarrow{ \sim } \ l \circ \id_A \
    \xRightarrow{ \id_l \times u } \ \, l \circ r \circ l \ \,
    \xRightarrow{ c \times \id_l } \ \, \id_B \circ \ l \, 
    \ \xRightarrow{ \sim } \ l \ )
    &&= \quad \id_l, \ \tag{zig} \label{eq:zig}
    \\
    & ( \ r \, \xRightarrow{ \sim } \ \id_A \circ \, r \
     \xRightarrow{ u \times \id_r } \  \, r \circ l \circ r \
    \xRightarrow{ \id_r \times c } \ r \circ \id_B 
    \ \xRightarrow{ \sim } \, r \ ) 
    &&= \quad \id_r. \ \tag{zag} \label{eq:zag}
 \end{align}
 We denote such an adjunction as $( l \adj r, u, c)$ and
 call $l$ the \emph{left adjoint} of $r$ and $r$ the \emph{right adjoint} of $l$.
 We say that 
 $u$ is the \emph{unit} and $c$ the \emph{counit} witnessing the adjunction $l \adj r$.
\end{defi}

\begin{exam}
A special case of this definition is the well known notion of adjoint functors. To see this, consider 
the bicategory $\ca C=\Cat_2$ whose objects are categories,  1-morphisms are functors,
and $2$-morphisms are natural transformations. Then the definition of a 1-morphism having a left/right adjoint in $\ca C$ is precisely that of having a left/right adjoint in the classical sense. \exampleend
\end{exam}

\begin{exam}
In the Morita bicategory $\Alg{1}$, objects are algebras, and, for two algebras $A$ and $B$, a $1$-morphism from $A$ to~$B$ is an $(A,B)$-bimodule $M$. The $1$-morphism $M$ has a left adjoint if and only if $M$ is finitely presented and projective as an $A$-module. We review this with details and proofs in \cref{sec:sec2}. \exampleend
\end{exam}

\begin{exam}
	Let $k$ be a field.
	The bicategory $\LinCat$ has as objects finite abelian $k$-linear categories, right exact linear functors
	as 1-morphisms and natural transformations as 2-morphisms. 	
	A 1-morphism $F$ has a left adjoint if and only if $F$ is left exact.	
	While every 1-morphism $G$ has a right adjoint linear functor $H$ it only is a right adjoint in $\LinCat$ if $G$ is right exact.
	For more details and proofs see \cite{douglas2018dualizable}. \exampleend
\end{exam}

For a bicategory $\ca C$ there is not just one opposite category but there are three different opposites.
The $1$-cell dual $\opcat{\ca C}$ reverses the direction of the $1$-morphisms but not the $2$-morphisms.
The $2$-cell dual $\cocat{\ca C}$ reverses the direction of the $2$-morphisms but not the $1$-morphisms.
The bidual $\coopcat{\ca C}$ reverses both $1$- and $2$-morphisms. By reversing the specified morphisms in \cref{Def:Adjunction} we can see what happens to an adjunction if we pass to one of these dual categories.

\begin{prop} \label{Prop:oppositeAdjunctionsBicat}
	Let $\ca C$ be a bicategory and $( l \adj r, u, c)$ be an adjunction in $\ca C$. We get the corresponding adjunctions
	\begin{displaymath}
		(r \adj l, u, c) \ \text{in} \ \opcat{\ca C}, \quad
		(r \adj l, c, u) \ \text{in} \ \cocat{\ca C} \quad
		\text{and} \quad
		(l \adj r, c, u) \ \text{in} \ \coopcat{\ca C}.
	\end{displaymath}
\end{prop}

A special case of adjunctibility are so called ambijunctions or ambidextrous adjunctions, where
the morphisms are simultanous left and right adjunctions.
These were first considered by Morita in \cite{Morita65}. A good reference is \cite{Lauda2005}.

\begin{defi}
    Let $\ca C$ be a bicategory. An \emph{ambidextrous adjunction} for a $1$-morphism $g$ are two adjunctions $f \adj g \adj h$
    with units and counits
    such that $f$ and $h$ are $2$-isomorphic, that is, there is an invertible 2-morphism between them.
\end{defi}

It is well-known that left and right adjoints of functors are unique in a certain sense. This carries over to our general setting with essentially the same proof, which we make precise below.
\begin{lema} \label{Lemma:adjointsUnique}
  Let $\ca C$ be a bicategory and $( f \adj h, u, c)$ and $( g \adj h, u', c')$ be adjunctions.
  Then the two left adjoints $f$ and $g$ of $h$ are $2$-isomorphic in $\ca C$. Moreover, the witnessing 2-isomorphism is unique, and the category of adjunctibility data is contractible.
  The same holds for right adjoints.
\end{lema}

Hence, having an adjoint is a property, although at first sight one must specify the data of the adjoint, unit, and counit.
\begin{defi}
  A $1$-morphism $f$ in $\ca C$ is \emph{left adjunctible} if it has a left adjoint $f^L$, i.e.~there exists an adjunction
  $f^L \adj f$.
  A $1$-morphism $f$ in $\ca C$ is \emph{right adjunctible} if it has a right adjoint $f^R$, i.e.~there exists an adjunction
  $f \adj f^R$.
  A $1$-morphism $f$ in $\ca C$ is \emph{adjunctible} if it is \emph{both} left and right adjunctible.
\end{defi}

In the rest of this subsection we record some Lemmas about properties of adjoints. Although well-known to experts, we include proofs of the following three lemmas in \cref{sec:appAdjoints} for reference .
\begin{lema} \label{Lemma:IsomorphismsOfAdjointsAreAdjunctible}
  Let $\ca C$ be a bicategory. Let $(l \adj r, u, c)$ be an adjunction in $\ca C$
  and $\mu: l \stackrel{\sim}{\Rightarrow} l'$ and $\nu: r \stackrel{\sim}{\Rightarrow} r'$ be $2$-isomorphisms.
  Then $(l' \adj r', u', c')$ is an adjunction with
  \begin{align}
    \label{eq:unitAndIso}
    \text{unit} \quad u' \ = \ 
    \big( \quad &\id \ 
    \xRightarrow{ \ u \ } \ 
    r \circ l \
    \xRightarrow{ \ \nu \times \mu \ } \ 
    r' \circ l' \ \big),
    \\
    \label{eq:counitAndIso}
    \text{counit} \quad c' \ = \ 
    \big( \quad & l' \circ r' \ 
    \xRightarrow{ \ \mu^{-1} \times \nu^{-1} \ } \ 
    l \circ r \
    \xRightarrow{ \ c \ } \ \id \ \big).
  \end{align}
\end{lema}

Note that we do not have uniqueness for the unit and counit in an adjunction.
\cref{Lemma:IsomorphismsOfAdjointsAreAdjunctible} shows that we can always change the units and counits by composition with an automorphism of the right (or left) adjoint to get another adjunction. This is in fact the only possibility.
\begin{lema} \label{Lemma:DifferentUnitsOfSameAdjunction}
    Let $\ca C$ be a bicategory. Let $(l \adj r, u, c)$ and $(l \adj r, u', c')$ be two adjunctions
    for the same $1$-morphisms $l, r$ in $\ca C$. Then there is a $2$-automorphism $\varphi:r \to r$, such that
    \begin{displaymath}
        (\varphi^{-1} \times \id_l) \circ u = u' 
        \qquad \text{and} \qquad
         c \circ (\id_l \times \varphi)  = c'.
    \end{displaymath}
\end{lema}

Next we record the fact that we can compose adjunctions and the resulting units and counits.
\begin{lema} \label{Lemma:composeAdjunctions}
    Let $\ca C$ be a bicategory. Let $f: X \to Y$ and $g: Y \to Z$
    be $1$-morphisms in $\ca C$ with adjunctions
    \begin{displaymath}
    ( f^L \adj f, u_f, c_f) \qquad \text{and} \qquad (g^L \adj g, u_g, c_g).
    \end{displaymath}
    Then 
        \begin{displaymath}
    ( f^L \circ g^L \adj  g \circ f, \mathtt{U}, \mathtt{C} ) \,,
    \end{displaymath}
where
    \begin{align*}
    \text{unit} \quad \mathtt{U} \ = \ 
    \big( \quad &\id_Z \ 
    \xRightarrow{ \ u_g \ } \ 
    g \circ g^L \
    \xRightarrow{ \ \id \times u_f \times \id \ } \ 
    g \circ f \circ f^L \circ g^L \ \big),
    \\
    \text{counit} \quad \mathtt{C} \ = \ 
    \big( \quad & f^L \circ g^L \circ g \circ f \ 
    \xRightarrow{ \ \id \times c_g \times \id \ } \ 
    f^L \circ f \
    \xRightarrow{ \ c_f \ } \ \id_X \ \big)
    \end{align*}
is an adjunction.

\end{lema}

\subsection{Adjoints in higher categories}
The generalization to adjoints of higher morphisms in higher categories is done by looking at their image in the relevant homotopy bicategories. For an informal explanation of homotopy bicategories see \cite[Remark 1.4.10]{lurie2009classification}; for a formal definition see \cite[Section 2.2.2]{CalaqueScheimbauer}, \cite{romo}.
\begin{defi}
  Let $\ca C$ be an $\infN N$category.
  A $1$-mor\-phism $f: X \to Y$ is \emph{left adjunctible}
  if $f$ is left adjunctible in the homotopy bicategory $\ho_2 \ca C$.
  For $k \geq 2$ fix two parallel $(k-2)$-morphisms $a$ and $b$. Let $V$ and $W$ be two
  $(k-1)$-morphisms from $a$ to $b$.
  A $k$-morphism $f: V \to W$ is \emph{left adjunctible}
  if $f$ is left adjunctible in the homotopy bicategory $\ho_2 \ca C(a,b)$, where
  $\ca C(a,b)$ is the $\Hom$-$\infN {N-k+1}$category from $a$ to $b$.
  The same definition holds for \emph{right adjunctible} with every
  ``left'' changed to~``right''.
\end{defi}

\begin{rema}
  Instead of using the homotopy categories we will work within the $\infN N$category
  $\ca C$ and check that for adjoints and duals the zig-zag identities
  \eqref{eq:zig} and \eqref{eq:zag} hold up to
  equivalence. This implies that they hold as equality in the corresponding homotopy bicategories.
\end{rema}

In an $\infN N$category we can extend the notions of
duals and adjoints to higher categorical versions. These 
do not only consist of one or a pair of adjunctions but instead
also demand that there are adjunctions for
the involved units and counits up to a certain level.

\begin{defi} \label{Def:nAdjunctibility} 
    Let $\ca C$ be an $\infN N$category and
    $f$ a $k$-morphism.
    A \emph{set of $1$-adjunctibility data} for $f$ is a left adjunction $(f^L \adj f, u_1, c_1)$ and
    right adjunction $(f \adj f^R, u_2, c_2)$ of $f$. 
    For $n\geq 2$ a \emph{set of $n$-adjunctibility data} is
    a set of $(n-1)$-adjunctibility data together with left and right adjunctions for all the unit and counit
    $(k + n - 1)$-morphisms in the $(n-1)$-adjunctibility data.
    A $k$-morphism $f$ is \emph{$n$-adjunctible} if there exists $n$-adjunctibility data for $f$.
    If all of the pairs of adjunctions in a set of $n$-adjunctibility data are ambidextrous
    we call it a \emph{set of ambidextrous $n$-adjunctibility data} und a morphism $f$ with such data
    \emph{ambidextrous $n$-adjunctible}.
\end{defi}

As is the case for adjunctibility, being $n$-adjunctible also is a property. This can be seen by showing that the $n$-adjunctibility data for a $k$-morphism $f$ is contractible, by using the following proposition.
The first part is stated as Lemma 4.1.5 in \cite{thesisAraujo}, which expresses the same proof in string diagrams.

\begin{prop}\label{Prop:uniqueness_of_adj}
    Let $\ca C$ be an $\infN N$category.
    Let $(l \adj r, u, c)$ and $(l' \adj r', u', c')$ be adjunctions of $k$-morphisms, such that $l \cong l'$ and
    $r \cong r'$ are isomorphic and $u$ is left adjunctible. Then $u'$ is also left adjunctible.
    The same is true for the counit and both statements also hold with  ``right'' instead of  ``left''.
    
    Hence, the choice of $(n-1)$-adjunctibility data of a morphism does not influence whether it is $n$-adjunctible or not.
\end{prop}
\begin{proof}
    Applying \cref{Lemma:IsomorphismsOfAdjointsAreAdjunctible} to both isomorphisms $l \cong l'$ and
    $r \cong r'$ we get an adjunction $(l \adj r, u'', c'')$, where $u''$ and $c''$ are compositions
    of isomorphisms with $u$ and $c$ as in \eqref{eq:unitAndIso} and \eqref{eq:counitAndIso}.
    By \cref{Lemma:DifferentUnitsOfSameAdjunction} we see that $u''$ and therefore also $u'$ are
    given by $u$ composed with an isomorphism.
    Hence, if $u$ is left adjunctible so is $u'$ by \cref{Lemma:composeAdjunctions}.
    
    The last part follows by induction.
\end{proof}

In $n$-adjunctibility data of a $k$-morphism $f$ we always have ``adjoint triples'' or ``two-step towers'' of adjunctions 

\[
\begin{array} {|r|rcl|}
	\hline
	k\text{-morphisms} &(f^L \adj f, u_1, c_1), &&(f \adj f^R, u_2, c_2)  \\ \hline
	(k+1)\text{-morphisms}  &u_1^L \adj u_1 \adj u_1^R,  &&c_1^L \adj c_1 \adj c_1^R, \\ 
	&u_2^L \adj u_2 \adj u_2^R,  &&c_2^L \adj c_2 \adj c_2^R \\ \hline
	\multicolumn{1}{|c|}{\vdots}
	&&\vdots &
	\\ \hline
	(k+n-1)\text{-morphisms} &
	\multicolumn{3}{c|}{ 2^{2n-1} \text{ adjunctions for all the units} } 
	\\ 
	&\multicolumn{3}{c|}{ \text{and counits from step $n-1$.} }   \\
	\hline
\end{array}
\]
which leads
to $2^{2(j-k) + 1}$ adjunctions of $j$-morphisms and $\frac{2}{3} \cdot (4^n-1)$
adjunctions overall. 

\subsection{Duals}

Finally, we define duals of objects and higher dualizability as a special case of the notions we have introduced so far. 

Given a monoidal category $\ca C$, we can construct the following bicategory $\mathbf{B}\ca C$, which is a one-object delooping. It has a unique object and the category of endomorphisms of this object is $\ca C$.  Composition of 1-morphisms in $\mathbf{B}\ca C$ is given by the monoidal product in $\ca C$.
In this setting, the notion of left and right adjoints in $\mathbf{B}\ca C$ can be translated to $\ca C$ as the notion of left and
right duals of objects, where the dual object is given by the
corresponding adjoint. We will work with symmetric monoidal categories
where the braiding leads to an isomorphism between left and right duals - see \cref{Lemma:DeloopingDualizableObject} - and hence, we do not have to distinguish between left and right.
Unravelling this, we arrive at the following definition.
\begin{defi}
  Let $(\ca C, \te, \uobj)$ be a symmetric monoidal category.
  An object $X \in \ca C$ is \emph{dualizable} or $1$-\emph{dualizable} if there is a \emph{dual} object
  $X^{*} \in \ca C$ and morphisms
  \begin{displaymath}
   u: \uobj \to X \te X^{*}, \qquad c: X^{*} \te X \to \uobj,
 \end{displaymath}
in $\ca C$ called \emph{unit} and \emph{counit}, satisfying the following identities, called \emph{zig-zag identities}
\begin{align}
    &(  X \phantom{^{*}}
    \cong \uobj \te X \phantom{^{*}}
     \xrightarrow{ \ u \te \id_X \ } \, \ X \te X^{*} \te X \,
    \xrightarrow{ \ \id_X \te c \  \, } \,  X \te \uobj  \phantom{^{*}} \cong  X 
     ) \ \, = \, \id_X,  \tag{zig} \label{eq:zigDual} \\
    &( X^{*} \cong X^* \te \uobj 
    \xrightarrow{ \id_{X^{*}} \te u }   \ X^{*} \te X \te X^{*} 
    \xrightarrow{ \ c \te \id_{X^{*}} \ }  \uobj \te X^{*} \cong  X^*
      ) = \, \id_{X^{*}}. \tag{zag} \label{eq:zagDual}
 \end{align}
\end{defi}
Using these notions we define higher dualizability for objects.
\begin{defi}  \label{Def:higherDual}
  Let $\ca C$ be a symmetric monoidal $\infN N$category.
  An object $X \in \ca C$ is \emph{$1$-dualizable} if it is so in the 
  homotopy category $\ho \ca C$.
  For $n \geq 2$ an object $X$ in $\ca C$ is \emph{$n$-dualizable} if it is 
  $1$-dualizable and the unit and counit in $\ca C$ are $(n-1)$-adjunctible.
\end{defi}

\begin{rema}
  The most well known notion of higher dualizability is that of
  fully dualizable objects as defined in \cite{lurie2009classification}.
  We will review it in \cref{sec:appDual} and see that for $n=N$ it
  is equivalent to \cref{Def:higherDual}.
\end{rema}

\section{Variants of adjunctibility} \label{sec:VariantsDual}

So far, for higher adjunctibility, we required the existence of both left and right adjoints of the morphism itself and all unit and counit morphisms. Now we vary this condition and allow for different choices. For instance, we could require that all morphisms  in question have only left adjoints; or have only right adjoints; or we could even demand infinitely long towers of adjunctions. These choices give rise
to a plethora of a priori different notions of higher adjunctibility. Our main result shows that this plethora in fact reduces to two cases.

 Of course, we have a good motivation for investigating these choices: the former two arise naturally when studying fully extended twisted/relative topological field theories \cite{JF_CS}, which we will recall in \cref{sec:applications}.
 This was our main motivation for this article.

As is the case for ordinary adjunctibility our
definitions will focus on the data involved
and define properties by demanding that such data exists.

\subsection{Higher left and right adjunctibility}\label{Sec:LR_adjunctibility}
We first recall the definitions of $n$-times left/right adjunctibility from \cite{JF_CS}.

\begin{defi}[\cite{JF_CS}, p.~194, before Theorem 7.6] \label{Def:higherAdj} 
    Let $\ca C$ be an $\infN N$category and
    $f$ a $k$-morphism.
    A \emph{set of $1$-times left adjunctibility data} for $f$ is
    a left adjoint $f^L$ together with unit $u$ and counit $c$ witnessing the adjunction between $f$ and $f^L$.
    For $n\geq 2$ a \emph{set of $n$-times left adjunctibility data} is
    a set of $(n-1)$-times left adjunctibility data together with left adjunctions for all the unit and counit
    $(k + n - 1)$-morphisms in the $(n-1)$-times left adjunctibility data.
    A $k$-morphism $f$ is \emph{$n$-times left adjunctible} if there exists
    a set of $n$-times left adjunctibility data.
    
    The same definition is made for \emph{$n$-times right adjunctibility} with every
    ``left'' changed to ``right''.
\end{defi}
Note that as for $n$-adjunctibility data, the category of such is contractible and hence, being $n$-times left adjunctible is a property.

Explicitly, a set of $n$-times left adjunctibility data for a $k$-morphism $f$ consists of adjunctions of
\[
\begin{array} {|r|ccc|}
	\hline
	k\text{-morphisms} &&&(f^L \adj f, u, c),  \\ \hline
	(k+1)\text{-morphisms}  &&&(u^L \adj u, u_u, c_u), \
	(c^L \adj c, u_c, c_c), \\ \hline
	(k+2)\text{-morphisms}  &&&(u_u^L \adj u_u, u_{uu},  c_{uu}), \
	(u_c^L \adj u_c, u_{uc},  c_{uc}), \\ 
	                    &&&(c_u^L \adj c_u, u_{cu},  c_{cu}), \
	(c_c^L \adj c_c, u_{cc},  c_{cc}). \\ \hline
	\vdots\hspace{1cm}
	&&&\vdots 
	\\ \hline
	(k+n-1)\text{-morphisms}  &&& 2^{n} 
	\text{ adjunctions for all the units}
	\\ 
	 &&&\text{and counits from step $n-1$.} \\
	 \hline
\end{array}
\]
For $k \leq j < k + n$ the number of adjunctions of $j$-morphisms is $2^{j-k}$. Thus $n$-times left adjunctibility data consists of $2^{n} -1$ adjunctions overall.

For examples we ask the reader to jump straight to  \cref{sec:sec2} where we take a look at $3$-times left and $3$-times right adjunctibility in a higher Morita category, which is a generalization of $\Alg1$.

\begin{rema}
    Note that $n$-times left together with $n$-times right adjunctibility does 
    not in general imply the stronger property of $n$-adjunctibility, since the
    data of left and right adjunctions do not need to be compatible.
    We will discuss this in more detail in  \cref{Theo:fromMixedToNadjunctibility}.
\end{rema}

Our main goal is to discuss when $n$-times left and $n$-times
right adjunctibility are equivalent. This is motivated by the following.
While ordinary left and right adjunctibility are in general vastly different notions we will see that this
difference vanishes for $2$-times left and $2$-times right adjunctibility.
The categorical fact responsible for this is the following Lemma, 
which can be found as Remark 3.4.22 in
\cite{lurie2009classification}. A proof is also given in
\cite[Lemma 1.4.4]{douglas2018dualizable}. 
We will give a slightly different proof, and be very explicit in exhibiting the necessary higher morphisms.

\begin{lema} \label{Lemma:MainLemma}(Interchange Lemma)
   Let $\ca C$ be an $\infN N$category.   Let $f$ be a $k$-morphism that admits 
    a left adjoint $f^L$ with unit $u$ and counit $c$. If $u$ and $c$ admit
    left adjoints $u^L$ and $c^L$, then 
    $(f \adj f^L, c^L, u^L)$ is an adjunction. Similarly, if $u$ and $c$ admit
    right adjoints $u^R$ and $c^R$, then 
    $(f \adj f^L, c^R, u^R)$ is an adjunction.
\end{lema}

 Before proving the Lemma, we first look at an immediate consequence.
\begin{coro}\label{cor:classsical_interchange}
Let $f$ be a $1$-morphism in an $\infN 3$category $\ca C$. Then $f$ is $2$-times left adjunctible if and only if $f$ is $2$-times right adjunctible.
\end{coro}

\begin{proof}
Applying the above Lemma to a set of $2$-times left adjunctibility data for $f$,
\begin{align}
&(f^L \adj f, u, c), \nonumber \\
&(u^L \adj u, u_u, c_u), \qquad
(c^L \adj c, u_c, c_c), \label{Eq:ucOf2LeftData}
\end{align}
we obtain the adjunction $(f \adj f^L, c^L, u^L)$, which
together with \eqref{Eq:ucOf2LeftData} forms $2$-times right adjunctibility
data for $f$.

Conversely, by \cref{Lemma:MainLemma},
we can also turn $2$-times right adjunctibility data
into $2$-times left adjunctibility data.
\end{proof}

Thus, $2$-times left and $2$-times right adjunctibility agree in contrast
to the classical notions of left- and right-adjunctibility.
Additionally, the adjunctions for $f$ above also form an ambidextrous adjunction, i.e.~a \emph{$2$-step tower}
of adjunctions $f^L \adj f \adj f^L$ for $f$ with the property
that the units and counits of the two adjunctions are adjoint to
each other. Of course, by symmetry we could continue this tower to a tower of infinite length
$$\cdots \adj f^L \adj f \adj f^L \adj f \adj f^L \adj \cdots$$

\begin{proof}[Proof of \cref{Lemma:MainLemma}]
    The adjunction $(f^L \adj f, u,c)$ consists of the data
    \begin{displaymath}
   		Y \stackrel[f]{f^L}{\RLarrows} X, \qquad u: \id_Y \Rightarrow f \  f^L,
    	\quad c: f^L  f \Rightarrow \id_X,
    \end{displaymath}
    with \eqref{eq:zig} and \eqref{eq:zag}
    \begin{align}
	    \label{Eq:3Lemma1}
	    \alpha: \ \big( \quad &f \, \ \ \xRightarrow{ \ u \times \id_f \ } \  \ f \ f^L f \
	    \xRightarrow{ \ \id_f \times c \ } \ f \quad \ \big)
	    &&\xRrightarrow{ \ \ \sim \ \ } \quad \id_f, \\
	    \label{Eq:3Lemma2}
	    \beta: \ \big( \quad &f^L \ \xRightarrow{ \id_{f^L} \times u } \ f^L  f \ f^L 
	    \xRightarrow{ c \times \id_{f^L} } \ f^L \quad  \big)
	    &&\xRrightarrow{ \ \ \sim \ \ } \quad \id_{f^L},
    \end{align}
    where $\alpha$ and $\beta$ are $(k+2)$-equivalences.
    
    To show that $(f \adj f^L, c^L, u^L)$ is an adjunction we need $(k+2)$-equivalences 
    \begin{align}
	    \label{Eq:3Lemma7}
	    \big( \quad &f^L \ \xRightarrow{ c^L \times \id_{f^L} } \ f^L f \ f^L 
	    \xRightarrow{ \id_{f^L} \times u^L } \ f^L  \quad  \big)
	    &&\xRrightarrow{ \ \ \sim \ \ } \quad \id_{f^L}, \\
	    \label{Eq:3Lemma8}
	    \big( \quad &f \, \ \ \xRightarrow{ \ \id_{f} \times c^L \ } \ f \ f^L f \ \,
	    \xRightarrow{ \ u^L \times \id_{f} \ } \ f \quad \ \, \big)
	    &&\xRrightarrow{ \ \ \sim \ \ } \quad \id_f.
    \end{align}
    
    We construct an equivalence for \eqref{Eq:3Lemma7}.
    From the adjunction data $(u^L \adj u, \eta_u, \epsilon_u)$ for the left adjoint of $u$ we obtain an adjunction
    \begin{equation} \label{Eq:3Lemma9}
	    ( \ \id_{f^L} \times u^L \ \adj \ \id_{f^L} \times u \, , 
	    \ \id_{\id_{f^L}} \times \eta_u \, , 
	    \ \id_{\id_{f^L}} \times \epsilon_u \ ).   
    \end{equation}
 	Similarly, from the adjunction data $(c^L \adj c, \eta_c, \epsilon_c)$ for the left adjoint of $c$ we obtain an adjunction
    \begin{equation} \label{Eq:3Lemma10}
	    ( \ c^L \times \id_{f^L} \ \adj \ c \times \id_{f^L} \, ,
	    \ \eta_c \times \id_{\id_{f^L}} \, , 
	    \ \epsilon_c \times \id_{\id_{f^L}} \ ).
    \end{equation}
    Both \eqref{Eq:3Lemma9} and \eqref{Eq:3Lemma10} are adjunctions in
    $\Hom_{\ca C}(X,Y)$.
    Applying \cref{Lemma:composeAdjunctions} to
    \begin{displaymath}
	    f^L \circ \id_Y
	    \ \stackrel[ \id_{f^L} \times u^L]{ \id_{f^L} \times u}{\LongDoubleRLarrows} \
	    f^L \circ f \circ f^L
	    \ \stackrel[ c^L \times \id_{f^L}]{c \times \id_{f^L} }{\LongDoubleRLarrows} \
	    \id_X \circ f^L
    \end{displaymath}
    yields an adjunction
    \begin{equation} \label{Eq:3Lemma11}
	    \big( \ 
	    ( \, \id_{f^L} \times u^L \, ) \circ ( \, c^L \times \id_{f^L} \, )
	    \ \adj \
	    ( \, c \times \id_{f^L}  \, ) \circ ( \, \id_{f^L} \times u \, ), \, \mathtt{U}, \, \mathtt{C}  
	    \ \big)
    \end{equation}
    with $\mathtt{U}$ and $\mathtt{C}$ given by \cref{Lemma:composeAdjunctions}.
    The right adjoint in \eqref{Eq:3Lemma11} is equivalent
    to $\id_{f^L}$ by assumption, via $\beta$ given by \eqref{Eq:3Lemma1}.
    By \cref{Lemma:IsomorphismsOfAdjointsAreAdjunctible} with $\mu= \id$ and $\nu = \beta$ we have an adjunction
    \begin{equation} \label{Eq:3Lemma12}
	    \big( \ 
	    ( \, \id_{f^L} \times u^L \, ) \circ ( \, c^L \times \id_{f^L} \, )
	    \ \adj \
	    \id_{f^L} , \, \mathtt{U}', \, \mathtt{C}'  
	    \ \big), 
    \end{equation}
    where the unit and counit are given by the compositions from \cref{Lemma:IsomorphismsOfAdjointsAreAdjunctible}.
    By the uniqueness of adjoints as stated in \cref{Lemma:adjointsUnique}
    and the fact that identities are self-adjoint
    the left adjoint in \eqref{Eq:3Lemma12} is equivalent to $\id_{f^L}$. We choose such an equivalence. This is the desired $(k+2)$-equivalence for \eqref{Eq:3Lemma7}.
    
   The construction of the equivalence for \eqref{Eq:3Lemma8} is analogous.
\end{proof}

\subsection{Mixed adjunctibility and reducing conditions}\label{Sec:mixed_adjunctibility}

In this subsection we will prove our main Theorem, which shows that higher left and right adjunctibility agree in even dimensions as a consequence of \cref{Lemma:MainLemma}.
In fact, we will prove a more general version for mixed versions of higher adjunctibility. Indeed, we will have $2^n$ a priori different notions of mixed adjunctibility, and our main Theorem reduces these cases to just two. In even dimensions, higher left and right adjunctibility reduce to the same case, whereas in odd dimensions they do not.

We start with mixed adjunctibility. We will do this by keeping track of when we require left and when we require right adjoints in a dexterity function. 

\begin{defi} \label{Defi:mixedAdj}
    Let $\ca C$ be an $\infN N$category.
    Let 
    \begin{displaymath}
    a^n: \{ 1, 2, \ldots, n \} \to \{ L, R \}. 
    \end{displaymath}   
    Let $f$ be a $k$-morphism in $\ca C$.
    A \emph{set of $a^n$-adjunctibility data} for $f$ with \emph{dexterity function $a^n$} is defined inductively as
    \begin{itemize}[itemsep=2mm]
        \item step $1$: the data of a
        $\begin{cases}
        	\text{left adjoint, } & \text{if} \quad a^n(1) = L, \\
        	\text{right adjoint, } & \text{if} \quad a^n(1) = R,
        \end{cases}$ \\
    	together with unit and counit witnessing the adjunction.
        \item step $j= 2, \ldots, n$: for
        all units and counits from step $j-1$
        \\ \phantom{step $1$:} the data of \phantom{a}
        $\begin{cases}
        	\text{left adjoints, } & \text{if} \quad a^n(j) = L, \\
        	\text{right adjoints, } & \text{if} \quad a^n(j) = R,
        \end{cases}$ \\
   		together with units and counits witnessing the adjunctions.
    \end{itemize}
    A $k$-morphism is \emph{$a^n$-adjunctible} if there exists a set of 
    $a^n$-adjunctibility data for~$f$.
\end{defi}

\begin{exam}
Let $f$ be an $n$-adjunctible $k$-morphism $f$ and $a^n$ any dexterity function. Then $f$ is $a^n$-adjunctible.
	\exampleend
\end{exam}

\begin{exam} \label{exam:a-adj}
    Let
    \begin{displaymath}
    l^n, r^n: \{ 1, 2, \ldots, n \} \to \{ L, R \}, 
    \qquad l^n(i)= L, \quad r^n(i) = R,
    \end{displaymath}
    be the constant dexterity functions with values `L' and `R', respectively. 
    Then $l^n$-adjunctibility is $n$-times left adjunctibility and
    $r^n$-adjunctibility is $n$-times right adjunctibility.
    
    We note that \cref{cor:classsical_interchange} shows that a 1-morphism is $l^2$-adjunctible if and only if it is $r^2$-adjunctible. Generalizing this statement is one of our main goals.
   \exampleend 
\end{exam}

\begin{exam}
Another natural pair of dexterity functions is given as follows.
One is simply the constant `R' dexterity function 
 \begin{equation} \label{eq:evenN}
 	\even^n = r^n
 \end{equation}
from \cref{exam:a-adj}.
The second one is given by the function
\begin{equation} \label{eq:oddN}
	\odd^n: \{ 1, 2, \ldots, n \} \to \{ L, R \}, 
	\quad \odd^n(n)= L, \quad \odd^n(i) = R \text{ for } 1 \leq i \leq n-1. 
\end{equation}
 	The names $\even$ and $\odd$ refer to the parity of the number of `L's in the image, which will become important below.
	\exampleend
\end{exam}

These two dexterity functions are natural in the following sense. Our main theorem proves that the $2^n$ a priori different notions of $a^n$-adjunctibility from \cref{Defi:mixedAdj} reduce to just two, namely exactly the $\even$ and $\odd$ ones.
The distinction for a given dexterity function $a^n$ is given by it being a parity dexterity of $\even^n$ or $\odd^n$ in the sense
of the following definition.

\begin{defi}
	Let $a^n, b^n: \{ 1, 2, \ldots, n \} \rightrightarrows \{ L, R \}$ be two dexterity functions.
	If
	\begin{displaymath}
		\mid (a^n)^{-1}(L) \mid \, \equiv \, \mid (b^n)^{-1}(L) \mid \mod 2,
	\end{displaymath}
	we say $b^n$ is a \emph{parity dexterity (function) of $a^n$} and the pair $a^n, b^n$ is a \emph{parity pair} of dexterity functions.
	If
	\begin{displaymath}
		\mid (a^n)^{-1}(L) \mid \, \not\equiv \, \mid (b^n)^{-1}(L) \mid \mod 2,
	\end{displaymath}
	we call $b^n$ a \emph{nonparity dexterity (function) of $a^n$} and the pair $a^n, b^n$ a \emph{nonparity pair} of dexterity functions.
\end{defi}

Now we formulate our main theorem that makes the equivalences of different $a^n$-adjunctibilities precise.

\begin{theo} \label{Prop:AnReducesToTwoClasses}
    Let $\ca C$ be an $\infN N$category and $n \geq 1$. \newline
    Let $a^n, b^n: \{ 1, 2, \ldots, n \} \rightrightarrows \{ L, R \}$ be two dexterity functions.
    \begin{enumerate}
    \item  If $a^n, b^n$ are a parity pair
    then a $k$-morphism $f$ is $a^n$-adjunctible if and only if $f$ is $b^n$-adjunctible.
    \item If $a^n, b^n$ are a nonparity pair
    then a $k$-morphism $f$ with an $a^n(1)$-adjoint $g$ is $a^n$-adjunctible if and only if $g$ is $b^n$-adjunctible.
    \end{enumerate}
\end{theo}
    
We will prove this theorem later in this section. 
Before that we record some immediate consequences. The first follows straight from the definitions of $\even^n$ and $\odd^n$.
\begin{coro}\label{cor:parity_implies_even/odd}
	Let $\ca C$ be an $\infN N$category and $n \geq 1$. \newline
	Let $a^n$ be a dexterity  function such that $f$ is $a^n$-adjunctible. If 
	 \begin{displaymath}
		| (a^n)^{-1}(L) | \, \equiv \, 0 \mod 2,
	\end{displaymath}
	then f is $\even^n$-adjunctible.
	If 
	\begin{displaymath}
		| (a^n)^{-1}(L) | \, \equiv \, 1 \mod 2,
	\end{displaymath}
	then f is $\odd^n$-adjunctible.
\end{coro} 

\begin{rema}\label{rem:exactly_two_classes}
One can construct examples of morphisms which are $\even^n$-adjunctible, but not $\odd^n$-adjunctible and vice versa. We will do so in Morita categories for $n=1$ in \cref{ex:rightnotleft1}, $n=2$ in \cref{Example:Alg2}, and $n=3$ in \cref{ex:3-times-left-not-right}. From these examples it is straightforward how to generalize to examples for arbitrary $n$ in the higher Morita categories $\Alg n$ from \cite{thesisScheimbauer}. Hence, there are precisely two independent (mixed) adjunctibility conditions we can impose on a $k$-morphism, which are $\even^n$ and $\odd^n$. 
\end{rema}

A special case of \cref{Prop:AnReducesToTwoClasses} is one of our main goals, namely, the equivalence between $n$-times left and right adjunctibility if $n$ is even.
\begin{coro} \label{Coro:evenDimensions}
    Let $\ca C$ be an $\infN N$category and $f$ a $k$-morphism.
    Let $n\geq 2$ be even. Then $f$ is $n$-times left adjunctible
    if and only if $f$ is $n$-times right adjunctible.
\end{coro}

Before we   prove \cref{Prop:AnReducesToTwoClasses} let us state a theorem
that describes the connection of the two classes of mixed adjunctibility to $n$-adjunctibility.
A $k$-morphism $f$ that is both left and right adjunctible is adjunctible by definition.
This generalizes in the following way.

\begin{theo} \label{Theo:fromMixedToNadjunctibility}
	Let $\ca C$ be an $\infN N$category and $n \geq 2$. 
	Let $a^n: \{ 1, 2, \ldots, n \} \to \{ L, R \}$ be a dexterity function and $f$ an $a^n$-adjunctible $k$-morphism.
	\begin{enumerate}
	\item \label{Theo:fromMixedToNadjunctibility1}	Then $f$ is ambidextrous $(n-1)$-adjunctible.
	\item  
	If $b^n$ is a nonparity dexterity of $a^n$ and $f$ also is $b^n$-adjunctible, then $f$ is $n$-adjunctible.
	\end{enumerate}	
\end{theo}
We will prove both this theorem and \cref{Prop:AnReducesToTwoClasses} later in this section.

As a consequence of the second statement we  see that there is a direct relation between $n$-adjunctibility and $a^n$-adjunctibility. However, although $n$-times left and $n$-times right adjunctibility represented by the dexterity functions $r^n$ and $l^n$ seemed like natural choices of one-sided adjunctibility, we see that 
 $\even^n$- and $\odd^n$-adjunctibility is a better choice, since they represent the two adjunctibility choices and together imply $n$-adjunctibility.
\begin{coro} \label{Coro:NadjunctibleIFFevenOddAdjunctible}
	Let $\ca C$ be an $\infN N$category and $n \geq 1$.
	A $k$-morphism $f$ is $n$-adjunctible if and only if it is both
	$\even^n$- and $\odd^n$-adjunctible.
\end{coro}

The first part of the theorem allows us to determine the adjunctibility type of the adjoint of a morphism.
\begin{coro} \label{Prop:AdjointIsbAdjunctible}	
    Let $\ca C$ be an $\infN N$category and $n \geq 2$. \newline
    Let $a^n, b^n: \{ 1, 2, \ldots, n \} \rightrightarrows \{ L, R \}$ be a nonparity pair of dexterity functions. 
    Let $f$ be an $a^n$-adjunctible $k$-morphism.
    Then any adjoint $g$ of $f$ is  $b^n$-adjunctible.
\end{coro}
\begin{proof}
By \cref{Theo:fromMixedToNadjunctibility} (\ref{Theo:fromMixedToNadjunctibility1}) and since $n\geq2$ the adjunction is ambidextrous, so we have that $g \adj f \adj g$.
We can extend this data to $a^n$-adjunctibility data for $f$.
If $a^n(1) = L$, then we use the adjunction $g\adj f$ and see that the $a^n$-adjunctibility data for $f$ gives $c^n$-adjunctibility data for $g$, where 
\[c^n(i) = \begin{cases} R & i=1,\\
					a^n(i) & i\neq 1. \end{cases} \]
Then by \cref{Prop:AnReducesToTwoClasses} (1) $g$ also is $b^n$-adjunctible. The case $a^n(1) = R$ is similar.
\end{proof}

Next, we obtain a corollary which generalizes the (quite obvious) statement that if a morphism $f$ is left adjunctible with left adjunctible left adjoint $g$, then $g$ is adjunctible.
\begin{coro} \label{Coro:adjointsOfNadjunctible}
	Let $\ca C$ be an $\infN N$category and $n \geq 1$.
	Let $a^n$ be a dexterity function. Let $g \adj f$ be an adjunction of $k$-morphisms such that both $f$ and $g$ are $a^n$-adjunctible.
	\begin{enumerate}
		\item 
		If $a^n(1)=L$ then $g$ is $n$-adjunctible,
		\item 
		If $a^n(1)=R$ then $f$ is $n$-adjunctible.
	\end{enumerate}
\end{coro}
\begin{proof}
	The case $n=1$ is immediate. For $n \geq 2$ we apply \cref{Prop:AdjointIsbAdjunctible} and obtain that $g$ is $b^n$-adjunctible for aa nonparity dexterity of $a^n$.
	By \cref{Theo:fromMixedToNadjunctibility} we have that $g$ is $n$-adjunctible. The case $a^n(1)=R$ is similar.
\end{proof}

\begin{coro} \label{Prop:AdjointAlsoAdjunctible}
    Let $\ca C$ be an $\infN N$category and $n \geq 2$.
    Let $f$ be an $n$-adjunctible $k$-morphism.
    Then any adjoint $g$ of $f$ is also $n$-adjunctible.
\end{coro}
\begin{proof}
Since $f$ is $n$-adjunctible, it is $a^n$-adjunctible for any dexterity function, for instance for $\even^n$ and $\odd^n$. Then by \cref{Prop:AdjointIsbAdjunctible} $g$ is $\odd^n$ and $\even^n$-adjunctible, respectively.
Hence, by \cref{Theo:fromMixedToNadjunctibility} (2) $g$ is $n$-adjunctible.
\end{proof}

\begin{nota}
	We denote the function on $\{R,L\}$ switching direction by
\begin{displaymath}
- \colon \{R,L\}\longrightarrow \{R,L\}, \qquad - L = R, \qquad - R = L.
\end{displaymath}
\end{nota}

To prove our theorems we will need a generalized version of the Interchange \cref{Lemma:MainLemma}. By applying the Interchange Lemma to one ``layer'' of $a^n$-adjunctibility data, we switch both the direction of the adjunctions in this and the consecutive ``layer'' as follows.
\begin{lema} \label{Lemma:GeneralInterchangeLemma}
Let $\ca C$ be an $\infN N$category and $n \geq 2$.
Let $a^n: \{ 1, 2, \ldots, n \} \to \{ L, R \}$ be a dexterity function.
Let $0 < j < n$ and define
 \begin{displaymath}
    b^n: \{ 1, 2, \ldots, n \} \to \{ L, R \}, \qquad i \mapsto
    \begin{cases}
      -{a^n(i)} \quad &i = j \text{ or } i =j+1, \\
      \phantom{-}a^n(i) \quad &i \neq j \text{ and } i \neq j+1.
   \end{cases}    
 \end{displaymath} 
 Then a $k$-morphism $f$ is $a^n$-adjunctible if and only if it is $b^n$-adjunctible.
 
In this case there is compatible $a^n$- and $b^n$-adjunctibility data
in the sense that
the adjunctions for $(k+j)$-morphisms together form $2$-step towers
\begin{displaymath}
	h \ \adj \ g \ \adj \ h\, ,
\end{displaymath}
such that the unit of the left (right) adjunction is adjoint to the counit of the right (left) adjunction, while all the other
adjunctions are the same.
\end{lema}
\begin{proof}
Let $a^n$ be a dexterity function and $f$ be a $k$-morphism that is $a^n$-adjunctible.
Let $0<j<n$. There are four different cases:
\vspace*{-0.5\multicolsep}
\begin{multicols}{2}
    \begin{enumerate}[label=(\roman*),itemsep=1mm]
        \item $a^n(j) = L, \ a^n(j+1) = L,$ \label{Item:DexterityLemmaSituationLL}
        \item $a^n(j) = R, \ a^n(j+1) = R,$ \label{Item:DexterityLemmaSituationRR}
    \end{enumerate}   
    \columnbreak   
    \begin{enumerate}[label=(\roman*),itemsep=1mm]
        \setcounter{enumi}{2}
        \item $a^n(j) = L, \ a^n(j+1) = R,$ \label{Item:DexterityLemmaSituationLR}
        \item $a^n(j) = R, \ a^n(j+1) = L$. \label{Item:DexterityLemmaSituationRL}
    \end{enumerate}
\end{multicols}
\vspace*{-0.5\multicolsep}
We explain the situation in the first case, the others are similar. Assume we have that $a^n(j) = a^n(j+1) = L$.
Let $g$  be a $(k+ j  - 1)$-morphism in set of $a^n$-adjunctibility data for $f$.  We have adjunctions of the form
\begin{align}
  & (g^L \adj g, u, c), \\
  & (u^L \adj u, u_u, c_u), \label{eq:adjLL2}
  && (c^L \adj c, u_c, c_c).  
\end{align}
Applying \cref{Lemma:MainLemma} we can switch the adjunction to an adjunction
\begin{displaymath}
  (g \adj g^L, c^L, u^L),
\end{displaymath}
exhibiting $g$ as having a {\em right} adjoint. Then \eqref{eq:adjLL2} exhibits right adjoints of the new unit $c^L$ and counit $u^L$. Hence, we have exchanged
the left adjoints with right adjoints.

We note that the new data fits 
with the rest of the data of adjunctions of $i$-morphisms
for $i<( k+ j  - 1)$ and $i> k+ j $.
This relies on the fact that the units
and counits $u_u, u_c, c_u, c_c$ of the adjunctions of $(k+ j )$-morphisms are
the same before and after the exchange, which proves the desired
compatibility.

Doing this for all 
$(k+ j  -1)$-morphisms in the $a^n$-adjunctibility data creates compatible
$b^n$-adjunctibility data, where $b^n$ is defined as above.

For the converse direction, start with a $k$-morphism $f$ that is $b^n$-adjunctible.
If we again assume that $a^n$ satisfies \ref{Item:DexterityLemmaSituationLL}, then $b^n$ satisfies \ref{Item:DexterityLemmaSituationRR} and so by swapping the role of $a^n$ and $b^n$ in the first part we obtain that $f$ is $a^n$-adjunctible. 
(Note that in this argument the cases \ref{Item:DexterityLemmaSituationLL} and \ref{Item:DexterityLemmaSituationRR} interchange roles; and \ref{Item:DexterityLemmaSituationLR} and \ref{Item:DexterityLemmaSituationRL} interchange roles.)
\end{proof}

Equipped with \cref{Lemma:GeneralInterchangeLemma} one could give a
direct proof of \cref{Coro:evenDimensions} by applying it to every second ``layer''
and iteratively switching `$LL \leftrightarrow RR$' in the constant dexterity functions
$l^n$ and $r^n$ from \cref{exam:a-adj}. Alternatively, we can use
\begin{displaymath}
	\mid (l^n)^{-1}(L) \mid \, = n \, \equiv \, 0 = \, \mid (r^n)^{-1}(L) \mid \mod 2
	\qquad \text{iff $n$ is even}
\end{displaymath}
and \cref{Prop:AnReducesToTwoClasses} which we proof now.

\begin{proof}[Proof of \cref{Prop:AnReducesToTwoClasses}]
	We start with the first claim.
	In the case $n=1$ two dexterity functions that form a parity pair are the same function and the statement is immediate.
	Let now $n \geq 2$.
	By \cref{Lemma:GeneralInterchangeLemma} we can switch
	consecutive entries
	\begin{equation} \label{eq:switch}
		\text{`}LL\text{'} \leftrightarrow \text{`}RR\text{'} \qquad \text{and} \qquad \text{`}LR\text{'} \leftrightarrow \text{`}RL\text{'}
	\end{equation}
	of the dexterity
	functions to get equivalent notions of $a^n$-adjunctibility.
	Note that both of these do not change the parity of the number of `$L$'s and `$R$'s.
	Let $a^n, b^n$ be a parity pair of dexterity functions.
	Let $a^n_1 := a^n$.
	For $j = 1, 2, \ldots, (n-1)$, we define a new dexterity function $a^n_{j+1}$ as follows. If $a^n_j(j) = b^n(j)$, set $a^n_{j+1} = a^n_j$.
If $a^n_j(j) \neq b^n(j)$, apply the exchange from
	\eqref{eq:switch} to entry $j$ and $j+1$ of $a^n_j$ to obtain $a^n_{j+1}$. By \cref{Lemma:GeneralInterchangeLemma}	a $k$-morphism $f$ is $a^n_{j+1}$-adjunctible if and only if it is $a^n_j$-adjunctible, and hence $a^n$-adjunctible.
	
After $n-1$ steps we constructed $a^n_n$ that describes an to $a^n$ equivalent notion
	of mixed adjunctibility and has $a^n_n(j) = b^n(j)$ for $1\leq j \leq n-1$ by construction.
	Since we did not change the parity of the number of `$L$'s in this process
	we have
	\begin{displaymath}
		|(a^n_n)^{-1}(L)| \, \equiv \, 
		| (a^n)^{-1}(L) | \, \equiv \,
		| (b^n)^{-1}(L) | \mod 2
	\end{displaymath}
	and therefore also $a^n_n(n) = b^n(n)$.
	This shows the equivalence of $a^n$- and $b^n$-adjunctibility for a parity pair.
	
	For the second part of the theorem let $a^n, b^n$ be a nonparity pair of dexterity functions.
	Let $f$ be an $a^n$-adjunctible morphism with $a^n(1)$-adjoint $g$. The case $n=1$ just states that left adjoints have right adjoints as well as the converse.
	For $n \geq 2$ it follows from the uniqueness of adjoints from \cref{Lemma:adjointsUnique} that $a^n$-adjunctibility data for $f$ is also
	$c^n$-adjunctibility data for $g$, where $c^n$ is the dexterity function with
	\begin{displaymath}
		c^n(j) := a^n(j) \text{ for } j = 2, 3, \ldots, n \text{ and } c^n(1) :=  -{a^n(1)}.
	\end{displaymath}
	Since $a^n,b^n$ and $a^n, c^n$ are nonparity pairs it follows that $b^n, c^n$ are a parity pair of dexterity functions.
	The claim now follows from the first part of the Theorem.
\end{proof}

At last we prove the second main Theorem. Part of the proof is the following Lemma.

\begin{lema} \label{Lemma:BuildUnitCounitTowersFromAmbidexterity}
	Let $\ca C$ be an $\infN N$category. Let 
	\[
		g \, \adj \, f \, \adj \, g 
	\]
	be an ambidextrous adjunction given by $(g \adj f, u_1, c_1)$ and $(f \adj g, u_2, c_2)$,
	such that $u_1$ and $c_1$ are $1$-adjunctible. Then $u_2$ and $c_2$ are
	also $1$-adjunctible.
\end{lema}
\begin{proof}
	Since $u_1$ and $c_1$ are $1$-adjunctible there are adjunctions
	\[
		u_1^l \adj u_1 \adj u_1^R \qquad \text{and} \qquad c_1^l \adj c_1 \adj c_1^R.
	\]
	Applying \cref{Lemma:MainLemma} for the left and the right adjoints of $u$ and $c$ yields
	\[
		(f \adj g, c_1^L, u_1^L) \qquad \text{and} \qquad (f \adj g, c_1^R, u_1^R).
	\]
	By \cref{Prop:uniqueness_of_adj} the morphisms $u_2$ and $c_2$ are left and right adjunctible.
\end{proof}

\begin{proof}[Proof of \cref{Theo:fromMixedToNadjunctibility}]
	Let $f$ be a $k$-morphism and $a^{n}$ a dexterity function such that $f$ is $a^{n}$-adjunctible.
	Consider a set of $a^{n}$-adjunctibility data.
	Apply \cref{Lemma:GeneralInterchangeLemma} with $j=1$ to
	get ambidextrous adjunctions for the $(k+1)$-morphisms.
	In particular, we have a set of ambidextrous $1$-adjunctibility data for $f$
	and its $(k+1)$-morphisms all have one adjoint $(k+1)$-morphism in the $a^{n}$-adjunctibility data
	by the compatibility condition in the Lemma.
	
	Inductively, in step $l=2, \ldots, n-1$ 
	assume we are given a set of ambidextrous $(l-1)$-adjunctibility
	data such that the $(k+l)$-morphisms have adjoints or are itself part of the $a^n$-adjunctibility data.
	Then we can apply \cref{Lemma:GeneralInterchangeLemma} with $j=l$ to $a^n$ and get ambidextrous adjunctions 
	for its $(k+l)$-morphisms. 
	By the induction hypothesis this extends the ambidextrous $(l-1)$- to ambidextrous $l$-adjunctibility
	data. Furthermore by the compatibility statement of \cref{Lemma:GeneralInterchangeLemma} the $(l+k+1)$-morphisms in this data
	all have adjoints in or are itself part of the $(l+k+1)$-morphisms of the $a^n$-adjunctibility data.
	After step $l = n-1$ we have constructed the desired ambidextrous $(n-1)$-adjunctibility data.
	
	For the second claim let $f$ also be $b^n$-adjunctible for a nonparity dexterity $b^n$ of~$a^n$.
	Consider another dexterity function $c^n$ given by
	\[
		c^{n}(i) = a^{n}(i) \ \text{ for } \ 1 \leq i \leq n-1 \quad
		\text{ and } \quad c^{n}(n) = -a^{n}(n).
	\] 
	We see that $c^n$ is a nonparity dexterity for $a^n$ and therefore a parity dexterity for $b^n$.
	By \cref{Prop:AnReducesToTwoClasses} we have that $f$ is $c^n$-adjunctible.
	Thus, we can take the corresponding adjunctions of the ambidextrous $(n-1)$-adjunctibility data and extend it
	to $c^{n}$-adjunctibility data. But since $c^{n}$ agrees with $a^n$ up to level $n-1$ we can take the same
	data and extend it to $a^{n}$-adjunctibility data. It follows from $c^{n}(n) = -a^{n}(n)$ that we have both left and
	right adjoints for the unit and counit $(k+n+1)$-morphisms of the $a^n$- and $c^n$-data.
	This gives us half of the adjunctions of $(k+n+1)$-morphisms that we need for $n$-adjunctibility.
	The missing half are adjunctions for the units and counits witnessing the adjunctions of $(k+n)$-morphisms that 
	do not appear in the $a^n$- and $c^n$-data. We get those
	from \cref{Lemma:BuildUnitCounitTowersFromAmbidexterity}.
\end{proof}

\subsection{Adjoints in opposite categories}\label{sec:opposite}
Similar to \cref{Prop:oppositeAdjunctionsBicat} for one adjunction let us record what happens if we pass to an opposite category. We will use this in \cref{sec:applications} for relative field theories.

\begin{defi}
Let $\ca C$ be an $\infN N$category. An \emph{opposite function} is a function
\begin{displaymath}
	\op^N: \{1, 2, \ldots, N\} \to \{ \id, \op \} \,.
\end{displaymath}
From an opposite function $\op^N$ we can construct the corresponding dual $\op^N$-category $\ca C^{\op^N}$, which reverses the direction of all $j$-morphisms with $\op^N(j) = \op$
for all $j$ but not any other morphisms.
\end{defi}

\begin{exam}\label{ex:odd/even op}
	Two opposite functions we will use later are
	\[
	\even \op (j) \coloneqq 
		\begin{cases}
				\op     \quad &j \text{ even}, \\
				\id 	&\text{else},
		\end{cases}
	\qquad \quad \text{and} \quad \qquad
	\odd \op (j) \ \coloneqq 
		\begin{cases}
				\op 	\quad &j \text{ odd}, \\
				\id     &\text{else\,.}
		\end{cases}
	\]
	\exampleend
\end{exam}

We will use the map
\begin{displaymath}
	\delta: \{ 0, 1 \} \to \{ \id, \op \}, \qquad 0 \mapsto \id, \ \ 1 \mapsto \op,
\end{displaymath}
to define an opposite function corresponding to a pair of dexterity functions in the following proposition.

\begin{prop} \label{Prop:oppositeFunctionAdjunction}
	Let $\ca C$ be an $\infN N$category and $f$ a $k$-morphism in $\ca C$.
	Let $a^n, b^n$ be two dexterity functions.
	Consider the opposite function 
	$\op_{a^n,b^n,k}: \{1, 2, \ldots, N\} \to \{ \id, \op \}$  given by
	\begin{align*}
		\op_{a^n,b^n,k}(j+k) &= \delta \left( \, | (a^j)^{-1}(L) + (b^j)^{-1}(L) | \mod 2 \, \right)
	\end{align*}
	for $j = 0,1, \ldots, n$,
	where all the other unassigned values are $\id$. Here $a^j$ denotes the dexterity function $a^n$ restricted to
	$\{1, \ldots, j \}$ with $a^0$ having empty domain.
	
	The following are equivalent:
	\begin{enumerate}
		\item  $f$ is $a^n$-adjunctible in $\ca C$ 
		\item  $f$ is $b^n$-adjunctible when considered as a $k$-morphism in $\ca C^{\op_{a^n,b^n,k}}$.
		\item  $f$ is $b^n$-adjunctible when considered as a $k$-morphism in $\ca C^{-\op_{a^n,b^n,k}}$, \\
		where $-\op_{a^n,b^n,k}$ has the opposite entries of $\op_{a^n,b^n,k}$.
	\end{enumerate}	
\end{prop}
\begin{proof}
	Since adjunctions in $\ca C$ are defined via the corresponding homotopy bicategory this follows from \cref{Prop:oppositeAdjunctionsBicat}.
\end{proof}

We now unravel this complicated looking statement in some examples.

\begin{exam} \label{exam:OddOppositeFunction}
	Let $\ca C$ be an $\infN N$category and let $f$ be a $1$-morphism which is $n$-times right adjunctible in $\ca C$. Set $a^n=r^n$.
	\begin{enumerate}		
	\item  Setting $b^n = a^n = r^n$, we have that
		\begin{displaymath}
			\op_{a^n,b^n,1} \equiv \id 
			\quad \text{and} \quad
			-\op_{a^n,b^n,1} \equiv \op
		\end{displaymath}
		are the constant functions valued $\id$ and $\op$, respectively.
		\item Setting $b^n=\odd^n$, we see that  $f$ is $\odd^n$-adjunctible in $\ca C^{\op_{a^n,b^n,1}}$ and in $\ca C^{-\op_{a^n,b^n,1}}$ with 
		\begin{displaymath}
			\op_{a^n,b^n,1}(j) =\begin{cases}
				\op 	& j = n+1, \\
				\id & \text{else\,}.
			\end{cases}
		\end{displaymath}
		
		\item 	Choosing $b^n=l^n$, we see that $f$ is $n$-times left adjunctible in $\ca C^{\op_{a^n,b^n,1}}$ and in $\ca C^{-\op_{a^n,b^n,1}}$ with 
		\begin{align*}
			\op_{a^n,b^n,1}(j) &= \begin{cases}
				\op     \qquad  &  j \leq n + 1, \ j \text{ even}, \\
				\id 	& \text{ else\,}.
			\end{cases}
		\end{align*}
		If $n+1 = N$ these are $\op_{a^n,b^n,1} = \even \op$ and $-\op_{a^n,b^n,1} = \odd \op$ from \cref{ex:odd/even op}.
	\end{enumerate}
	\exampleend
\end{exam}

\section{Adjunctibility in Morita categories} \label{sec:sec2}

In this section we take a look at the Morita bicategory and similar Morita 
$3$- and $4$-categories. We construct examples of morphisms which satisfy different notions of adjunctibility introduced in the previous section.

In \cite{Morita1958DualityFM}
Morita introduced what is now called Morita equivalences between rings.
In \cite{IntroToBiCat} B\'{e}nabou recognized that these are equivalences in a bicategory, the Morita bicategory $\Alg 1$.
It can be described as follows. Let $k$ be a commutative ring or field.
\begin{itemize}
  \item objects are associative $k$-algebras,
  \item for two $k$-algebras $A$ and $B$ the $1$-morphisms $A \to B$ are 
  bimodules of the form $_A M_B$. Composition of two bimodules $_A M_B$ and
  $_B N_C$ is given by the relative tensor product
  \begin{displaymath}
    {_B N_C  }\circ { _A M_B } =  {_A M_B } \ {\te_B } \ { _B N_C },
  \end{displaymath}
  which is an $(A,C)$-bimodule,
  \item $2$-morphisms are given by bimodule homomorphisms.
\end{itemize}

Recall that $_A M_B$ means that $M$ is a left $A$-module and right $B$-module
such that $(am)b = a(mb)$ holds for all $a \in A, b \in B$ and $m \in M$.
Alternatively, a bimodule $_A M_B$ is the same thing as a left module over the $k$-algebra
$A \te_k B^{\op}$, where $B^{\op}$ is the opposite algebra of $B$.
A bimodule homomorphism $_A M_B \to {_A N_B}$ is a linear map which is both a left $A$-module
and a right $B$-module homomorphism $M \to N$.
Together with the tensor product of $k$-algebras $\Alg 1$ is a symmetric
monoidal bicategory.
Full details can be found in \cite{IntroToBiCat}.

\begin{rema}
\emph{Warning:} Sometimes a different convention for the direction of 1-morphisms is chosen, namely, that $_A M_B$ is a 1-morphism from $B$ to $A$.
\end{rema}

To characterize $2$-dualizability we need to know which
$k$-algebras are $1$-dualizable and which bimodules have left and which have right adjoints.
These results are well-known, but we include proofs in \cref{sec:appMorita} for the interested reader.

\begin{prop} \label{Prop:allAlgDual}
  Every object in $\Alg 1$ is $1$-dualizable. For a $k$-algebra $A$ its dual is the opposite algebra $A^{\op}$,
  the unit morphism is the bimodule $_k A_{A \te_k A^{\op}}$ and the counit is
  $_{A \te_k A^{\op}} A _k$.
\end{prop}

\begin{prop} \label{Prop:bimodDual}
  A bimodule $_A R_B$ has a left adjoint if and only if $R$ is finitely presented
  projective as a left $A$-module. In this case a dual is $_B L _A$, where
  $L \cong \Hom_A(R, A)$.
  
  Dually,
  a bimodule $_B L_A$ has a right adjoint if and only if $L$ is finitely presented
  projective as a right $A$-module. In this case the dual is $_A R_B$, where
  $R \cong \Hom_A(L, A)$.
\end{prop}

\begin{exam}\label{ex:rightnotleft1}
With this proposition, we immediately find examples of right adjunctible but not left adjunctible 1-morphisms in $\Alg 1$:
Take any $k$-algebra $M$ which is a left $A$-module and as such is not finitely presented projective.
Then $M$ is an $(A,M)$-bimodule. Since $M$ as a trivial $M$-module is
finitely presented projective by \cref{Prop:bimodDual} we have that 
the morphism $_A M_M$ is right adjunctible but not left adjunctible.
Likewise, $_M M_{A^\op}$ is left adjunctible but not right adjunctible.
Note that in addition by \cref{Prop:allAlgDual} both $A$ and $M$ are $1$-dualizable.
This is important when interpreting this as a relative TFT.
	\exampleend
\end{exam}

There is a natural generalization of the Morita bicategory to a higher
category ${\Alg n} (\ca C)$ whose objects are $E_n$-algebras in some suitable symmetric monoidal (possibly higher) category $\ca C$. 
If $\ca C$ is a symmetric monoidal $(\infty,k)$-category satisfying some mild conditions, then ${\Alg n} (\ca C)$ is a symmetric monoidal $(\infty,n+k)$-category, sometimes called a Morita $(n+k)$-category. See \cite{Haugseng2017}, \cite{thesisScheimbauer}, \cite{JF_CS} for the technical details.
The classical Morita bicategory described above is the special case $n=1$ and $\ca C=\Vect_k$.

In the rest of this section we will look at two other special cases in which the construction reduces to familiar objects: In both cases we take $\ca C=\Vect_k$. We first look at the (weak) 3-category ${\Alg 2} = {\Alg 2}(\Vect_k)$ and then we consider the  (weak) $4$-category ${\Alg 3} = {\Alg 3}(\Vect_k)$.

Objects in ${\Alg 2}$ and ${\Alg 3}$ are simply commutative $k$-algebras. Certain morphisms are bimodules which themselves have an algebra structure. These satisfy a compatibility condition as follows:
\begin{defi}
Let $C$ be a commutative algebra and $M$ an algebra, which also has a $C$-module structure. Then the action factors through the unit of $M$,
\[C\otimes M \longrightarrow  M, \qquad c\otimes m \longmapsto c\cdot m = c\cdot 1_M \cdot m,\]
so the action is determined by a homomorphism
\[C\longrightarrow M, \qquad c \longmapsto c\cdot 1_M \,.\]
We say that  $M$ is a {\em central $C$-module} if this map factors through the center $Z(M)$.
An algebra which is an $(A,B)$-bimodule is central if it is central as an $A\otimes B^{\op}$-module.
\end{defi}

\begin{exam}
If $M$ is commutative, any module is central.
	\exampleend
\end{exam}

We start describing the weak 4-category $\Alg 3$ and state some dualizability and adjunctibility results. Rather than repeating the same statements and arguments for $\Alg2$, we use that $\Alg 2 = \Hom_{\Alg3} (k,k)$ afterwards. The reader feeling more comfortable starting with a lower-dimensional setting can jump to \cref{Example:Alg2} first and work their way backwards.

Informally, $\Alg 3$ can be described as follows.
\begin{itemize}
  \item Objects are commutative $k$-algebras,
  \item a $1$-morphism between two commutative algebras $T \to S$ is a
   commutative $k$-algebra $A$ which also is a
    bimodule $_T A_S$;
   \item a $2$-morphism between two $1$-morphisms $A \to B$ is an associative $k$-algebra $M$ which also is a
   central bimodule $_A M_B$;
   \item a $3$-morphism between two $2$-morphisms $N \to M$ is a bimodule $_M \Sigma_N$ such that for every $\sigma\in \Sigma$ and $c\in A\otimes B^{\op}$, we have that 
   \begin{equation}\label{Eqn:centeringBimodAlg3}
   \left( c \cdot 1_M\right)\cdot \sigma = \sigma \cdot \left(c \cdot 1_N\right)\,;
   \end{equation}
   \item and a $4$-morphism is a homomorphism, meaning it is compatible with all actions of the source and target objects, 1-morphisms, and 2-morphisms.
\end{itemize}
Compositions of morphisms which are bimodules is given by the relative tensor product, similarly to the composition of 1-morphisms in $\Alg 1$. Composition of 4-morphisms is composition of homomorphisms.
\begin{rema}
The structure of a 3-morphism as defined above automatically is compatible with the action of the source and target objects, meaning that the actions on $\Sigma$ also commute as in \eqref{Eqn:centeringBimodAlg3}.
\end{rema}
\begin{nota}
We iterate the subscript notation for bimodules as a compact notation for higher morphisms. For instance, for a commutative $k$-algebra $T$ the identity
$3$-morphism on the identity of the identity of $T$ is
\begin{displaymath}
  \id_{\id_{\id_T}} \ = \  _{ _{_T T_T } T_{ _T T_T } } T_{ _{ _T T_T } T_{ _T T_T } }.
\end{displaymath}
\end{nota}

We will need two basic facts about existence of duals and adjoints in $\Alg3$, described in the following two propositions.
The first is a generalization of \cref{Prop:allAlgDual}.
\begin{prop}[Special case of Theorem 5.1 in \cite{GS}]\label{Prop:Alg3dualizable}
The symmetric monoidal weak 4-category $\Alg3$ is fully 3-dualizable, i.e.
\begin{enumerate}
\item  every object $T$ in $\Alg3$ is 3-dualizable with dual $T^{\op}=T$; and
\item any $1$-morphism $_TA_S$ and any $2$-morphism $_AM_B$ in $\Alg3$ has both a left and a right adjoint whose underlying algebras are $A^{\op}=A$ and $M^{\op}$.
\end{enumerate}
\end{prop}

\begin{proof}
We spell out the data in the case of a $1$-morphism $T \to S$ for later reference. The proof of the other statements is similar.

Let $T$ and $S$ be commutative $k$-algebras, and let $A$ be a commutative $k$-algebra which is a bimodule $_T A_S$.
By commutativity we have that $A^{\op} = A$ and hence $A$ also is
an $_S A_T$ bimodule.

A right adjoint of $_T A_S$ is given by $_S A_T$ with
\begin{align*}
  & \phantom{\text{co}}\text{unit} \qquad 
  ( \, u : \id_T \to { _T A_S } \, {\te_S} \, _S A_T \, ) \ 
  = \  _{ _T T_T } A_{ ( { _T A_S } ) { \, {\te_S} \, } ( { _S A_T } ) } , 
   \\ \\ & \text{counit}
  \qquad ( \, c : { _S A_T } \, {\te_T} \, _T A_S \to \, \id_S \, ) \ 
  = \  _{ ( { _S A_T } ) \, {\te_T} \, ( { _T A_S } ) } A_{ _S S_S }.
\end{align*}
For the unit $u$ we have an adjunction $(u \adj u^R, u_u, c_u)$
with
\begin{align}
  \label{eq:CoUnitsAlg3_1}
  & \phantom{\text{counit}} \qquad 
  u^R && = \
   _{  ( { _T A_S } ) { \, { \te_S } \, } ( { _S A_T } ) } A _{ _T T_T } , 
  \\ \nonumber \\ \label{eq:CoUnitsAlg3_2}
  & \phantom{\text{co}}\text{unit} \qquad 
    u_u && = \ 
    _{ T } A_{ { ( _T A_{A { \te_S } A} ) } 
    \ { \te_{( A { \te_S } A )} } \ 
    { ( _{A { \te_S } A} A_T ) }  }
    ,    
   \\ \nonumber \\ \label{eq:CoUnitsAlg3_3}
   & \text{counit} \qquad
   c_u && = \
   _{ { ( _{A { \te_S } A} A_T ) }  \ {\te_T} \ { ( _T A_{A { \te_S } A} ) } } 
   A_{ _{A {\te_S } A} ( A {\te_S } A)_{A {\te_S } A}  }
   .
\end{align}
For the counit $c$ we have an adjunction $(c \adj c^R, u_c, c_c)$
with
\begin{align}
  \label{eq:CoUnitsAlg3_4}
  & \phantom{\text{counit}} \qquad 
  c^R && = \
   _{ _S S_S } A _{  ( { _S A_T } ) \, {\te_T} \, ( { _T A_S } ) } , 
  \\ \nonumber \\ \label{eq:CoUnitsAlg3_5}
  & \phantom{\text{co}}\text{unit} \qquad 
    u_c && = \ 
    _{ { ( _S A_{A { \te_T } A} ) } 
    \ { \te_{ ( A { \te_T } A ) } } \ 
    { ( _{A { \te_T } A} A_S ) } }
    A_{ { ( _{ { A \te_S } A } A_S ) } \ { \te_S } \ { _S A_{ A {\te_T } A } } } 
    , 
    \\ \nonumber \\ \label{eq:CoUnitsAlg3_6}
    & \text{counit} \qquad
   c_c && = \
   _{ { ( _{A { \te_T } A} A_T ) }  \ { \te_T} \ { ( _T A_{A { \te_T } A} ) } } 
   A_{ S }
   .
\end{align}
Note that we shortened the notation of the $3$-morphisms 
$u_u, u_c, c_u, c_c$ and wrote $A$ instead of $_T A_S$ or $_S A_T$
on the lowest level. Similarly, the $T$ and $S$ in $u_u$ and $c_c$ are actually
\begin{displaymath}
  _{_T T_T} T_{_T T_T} \qquad \text{and} \qquad _{_S S_S} S_{_S S_S}.
\end{displaymath} 
All the adjunctions above form $2$-times right adjunctibility data for $_T A_S$.
\end{proof}

The second statement we need is a generalization of \cref{Prop:bimodDual}.
Let $T,S$ be objects, $A,B$ be 1-morphisms from $T$ to $S$, $M, N$ 2-morphisms from $A$ to $B$, and $_M \Sigma_N$ from $M$ to $N$ a 3-morphism in $\Alg3$.
\begin{prop}\label{Prop:bimodAdjointAlg3}
 The 3-morphism $_M \Sigma_N$ in $\Alg3$ has a left adjoint if and only if $\Sigma$ is finitely presented and projective merely as a module over the underlying algebra of $M$, and similarly for the right adjoint. The adjoints are given by $\Hom_M(\Sigma, M)$ and $\Hom_N(\Sigma, N)$, respectively.
\end{prop}

For this we need a statement about ``lifting'' the adjoints from \cref{Prop:bimodDual} to adjoints of 3-morphisms.
This lifting procedure is exactly the same as the one described in \cite[Proposition 5.17]{BJS}.

\begin{prop}\label{Prop:lifting}
The $(N,M)$-module $\Hom_M(\Sigma, M)$ has a canonical structure of an $(A,B)$-bimodule, given for $F\in \Hom_M(\Sigma, M)$ and $c\in A\otimes B^{\op}$ by
\[ c \cdot F (-) \coloneqq F ( - \cdot (1_Nc)) =  F ((1_Mc) \cdot -) = (1_Mc) \cdot  F ( -) = F ( -) \cdot (1_Mc) \eqqcolon  F( - ) \cdot c \,.\]
Similarly, this structure is compatible with the structure of $A,B$ being $(R,S)$-modules and hence $\Hom_M(\Sigma, M)$ is a 3-morphism from $N$ to $M$.
\end{prop}

\begin{proof}
The first equality holds by \eqref{Eqn:centeringBimodAlg3}. The second holds because $F$ is an $M$-module. The third equality holds because the action of $c$ on $1_M$ lands in the center~$Z(M)$.
\end{proof}

\begin{proof}[Proof of \cref{Prop:bimodAdjointAlg3}]
We discuss the case of the left adjoint; the statement for the right adjoint is similar.
Let $_M \Sigma_N$ be a 3-morphism in $\Alg 3$, where $M$ and $N$ are $(A,B)$-bimodules, and $A$ and $B$ are $(T,S)$-modules.
There is a natural map from the (homotopy) bicategory $\ho_2 \ca C(A,B)$ to $\Alg 1$, forgetting the extra bimodule structures.
We can apply \cref{Prop:bimodDual} to the image of $_M \Sigma_N$ and obtain a left adjoint $\Hom_M(\Sigma, M)$ in $\Alg 1$ if and only if $\Sigma$ is finitely presented and projective over $M$. Now applying \cref{Prop:lifting} $\Sigma$ can be given the structure of a 3-morphism from $N$ to $M$. The unit and counit of the adjunction in $\Alg 1$ are homomorphisms of the underlying bimodules, but since the actions of $T,S,A,B$ are all central, the unit and counit respect this action as well.
\end{proof}

\begin{exam}\label{ex:3-times-left-not-right}
Now we take $T,S$ commutative $k$-algebras such that $A=T$ is a right $S$-module which is not finitely presented and projective.
We claim that in $\Alg3$, as a 1-morphism
$_A A_S$ is 2-times right adjunctible, but not $3$-times right adjunctible. Moreover, it is 3-times left adjunctible.

To see this, in \cref{Prop:Alg3dualizable} we can read off the $2$-times right adjunctibility data of the adjunction
\begin{equation}\label{eq:adjBimod}
_AA_S \adj {_SA_A}\,.
\end{equation}
First look at the counit $c_c$ in \eqref{eq:CoUnitsAlg3_6}. By \cref{Prop:bimodDual} it does not have a right adjoint. 
So we cannot complete the $2$-times right adjunctibility data to
 $3$-times right adjunctibility data. By the uniqueness of adjoints up to isomorphism (see \cref{Prop:uniqueness_of_adj})
 it follows that $_A A_S$ is not $3$-times right adjunctible.
 
To look at left adjunctibility we could take the general adjunctibility data of
\begin{displaymath}
  ( { _T A_S} \adj {_S A_T}, u, c )
\end{displaymath}
from above and swap $S$ and $T$ to obtain 2-times left adjunctibility data of $_TA_S$.

Instead we apply \cref{Lemma:GeneralInterchangeLemma} to the adjunction
of $1$-morphisms \eqref{eq:adjBimod}
to obtain the adjunction
\begin{displaymath}
  ( { _S A_A} \adj {_A A_S } , c^R, u^R ).
\end{displaymath}
Left adjoints for $c^R$ and $u^R$ are given by \eqref{eq:CoUnitsAlg3_1} and
\eqref{eq:CoUnitsAlg3_4} with $T=A$.
Units and counits of these adjunctions 
are given by $u_u, u_c, c_u$ and $c_c$ with $T=A$ in 
\eqref{eq:CoUnitsAlg3_2}, \eqref{eq:CoUnitsAlg3_3}, \eqref{eq:CoUnitsAlg3_5}
and \eqref{eq:CoUnitsAlg3_6}.
These are all $A$, viewed as a left module for --  respectively --
\begin{displaymath}
A,
\qquad
A { \te_A } A,
\qquad
A { \te_{ ( A \ { \te_A } \ A  )} A }
\qquad \text{and} \qquad
A { \te_A } A\,.
\end{displaymath}
Since in each case it is finitely presented projective, $u_u, u_c, c_u$ and $c_c$ all have left adjoints by
\cref{Prop:bimodAdjointAlg3}. This is a $3$-times left adjunctibility datum.
Additionally, $A$ and $S$ are $3$-dualizable in $\ca C$, since all
objects of $\Alg 3$ are $3$-dualizable. \exampleend
\end{exam}

\begin{rema}
It is possible to generalize this idea and find examples in all higher categories ${\Alg n} (\Vect_k)$ for $n$ odd.
\end{rema}

Morita categories also allow us to look at an example of mixed adjunctibility in the sense of
\cref{Defi:mixedAdj}. For $n=2$ we have $4$ possible dexterity functions $a^n$ given by `$LL$', `$RR$'
`$RL$' and `$LR$', where the first two and the latter two define equivalent notions of adjunctibility by
\cref{Prop:AnReducesToTwoClasses}. To see that in general these two classes are different we
consider $1$-morphisms in the weak $3$-category $\Alg 2$.
We could define $\Alg 2 = \Hom_{\Alg3} (k,k)$, but let us unravel this.
Informally, $\Alg 2$ can be described as follows.
\begin{itemize}
  \item Objects are commutative $k$-algebras,
  \item a $1$-morphism between two $1$-morphisms $A \to B$ is an associative $k$-algebra $M$ which also is a
   central bimodule $_A M_B$;
   \item a $2$-morphism between two $2$-morphisms $N \to M$ is a bimodule $_M \Sigma_N$ such that for every $\sigma\in \Sigma$ and $c\in A\otimes B^{\op}$, we have that 
   \begin{equation}\label{Eqn:centeringBimodAlg2}
   \left( c \cdot 1_M\right)\cdot \sigma = \sigma \cdot \left(c \cdot 1_N\right)\,;
   \end{equation}
   \item and a $3$-morphism is a homomorphism, meaning it is compatible with all actions of the source and target objects and 1-morphisms.
\end{itemize}
Compositions of morphisms which are bimodules is given by tensoring over the middle
algebra, similarly to the composition of 1-morphisms in $\Alg 1$. Composition of 3-morphisms is composition of homomorphisms.

Let $T$ and $S$ be commutative $k$-algebras and $_T A_S$ be a $1$-morphism $T \to k$.
Similar to $\Alg 3$ we have an adjunction
$( {_T A_S} \adj {{}^{\vphantom{\op}}_S A^\op_T}, u, c)$ with
\begin{align} 
	& \phantom{\text{co}}\text{unit} \qquad 
	( \, u : \id_T \to { _T A_S } \, {\te_S} \, _S A_T \, ) \ 
	= \  _{ _T T_T } A_{ ( { _T A_S } ) { \, {\te_S} \, } ( { _S A_T } ) } , 
	\label{eq:unitAlg2}
	 \\ & \text{counit}
	\qquad ( \, c : { _S A_T } \, {\te_T} \, _T A_S \to \, \id_S \, ) \ 
	= \  _{ ( { _S A_T } ) \, {\te_T} \, ( { _T A_S } ) } A_{ _S S_S }.
	\label{eq:counitAlg2}
\end{align}

\begin{exam}\label{Example:Alg2}
	Consider the $1$-morphism $_T A_S$ in $\Alg 2$ with $A = S = k$ and $T = k[x] / (x^2)$ acting trivially.
	Recall that the dexterity function $\odd^2$ is given by $\odd^2(1) = R$ and $\odd^2(2) = L$.
	We have data of a right adjoint of $_T A_S$ as given above. The question is whether we can extend this to $\odd^2$-adjunctibility data. This is not the case,
	since $A$ is not projective as a left $T$-module and by \cref{Prop:bimodDual} the unit $u$ in
	\eqref{eq:unitAlg2} is not left adjunctible.
	Thus $_T A_S = {}_{k[x] / (x^2)} k_k$ is not $\odd^2$-adjunctible.
	
	On the other hand, $A = k$ is trivially a finitely presented projective module over $k = k_k \te_{k} {}_k k^{\op}$. Thus, by \cref{Prop:bimodDual} the unit $u$ in
	\eqref{eq:unitAlg2} and the counit $c$ in \eqref{eq:counitAlg2} with $A = S = k$ have right adjoints.
	Since $\even^2(1) = \even^2(2) = R$,
	we get that $_T A_S = {}_{k[x] / (x^2)} k_k$ is $\even^2$-adjunctible.	
	\exampleend
\end{exam}

\section{Applications to field theories} \label{sec:applications}

\newcommand{\oppp}{\mathrm{oppp}}

One of our motivations for this article was to understand the dualizability requirements for the notion of a relative version of topological field theories. In this section we explain our results in this context and prove some variations and extensions.

 Relative versions of field theories were first proposed by Stolz--Teichner in \cite{stolz2011supersymmetric} using the name {\em twisted field theory}, and by Freed--Teleman in \cite{FT2014relative} using the name {\em relative field theory}. More recently, Freed--Moore--Teleman \cite{FMT} introduced the term {\em quiche} as being half of a {\em sandwich field theory}, which in turn has a boundary on two sides; both in the context of symmetries of field theories, see also \cref{foot_Butterbrot}. 

Relative theories should capture the following idea. Instead of attaching a number to top-dimensional closed manifolds, we choose an element in a vector space (for instance a line), implemented as a morphism $\uobj \to V$, and this vector space should depend on the manifold $M$,
\begin{equation}\label{eqn:element_in_vsp}
\uobj \to T(M) \,.
\end{equation}
A first approximation to the notion of a relative field theory, motivated by physical bulk-boundary systems, is an \emph{(n+1)-dimensional field theory with boundaries}, that is, a symmetric monoidal functor
\[\TFT{Z}^\partial\colon \Bord_{n+1}^{\partial} \longrightarrow \mathcal{C},\]
where $\Bord_{n+1}^{\partial}$ is a (possibly higher) category\footnote{Such a higher category is outlined in \cite{FreedTelemanGapped, FMT} and was constructed in detail by William Stewart in his PhD thesis \cite{StewartThesis}.} of cobordisms with free (=marked) boundaries, and $\mathcal{C}$ is a symmetric monoidal $(\infty,N)$-category for $N>n$. In the fully extended framed situation, an extension of the Cobordism Hypothesis, namely Lurie's Cobordism Hypothesis with singularities \cite{lurie2009classification} explained that these are fully determined by their value at an interval with one incoming boundary component (black) and one free boundary (red) viewed as a bordism from one point to the empty set:
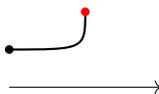
\[
	\begin{tikzpicture}
		\draw[thick] (0,0) ..controls (1,0).. (1,0.5);	
		\filldraw (0,0) circle (0.05cm);
		\filldraw[red] (1,0.5) circle (0.05cm);
		
		\draw[->] (0,-0.5) -- (2,-0.5);
	\end{tikzpicture}
\]
This value is an $n$-adjunctible 1-morphism from an $(n+1)$-dualizable object to the unit in $\mathcal C$. In turn, according to the Cobordism Hypothesis with singularities, every such $n$-adjunctible 1-morphism gives rise to a framed fully extended $(n+1)$-dimensional TFT with boundaries, and in particular a framed fully extended $(n+1)$-dimensional TFT $\TFT{Z}^{bulk}$ by restricting to bordisms with no free (red) boundary via the inclusion
\[\Bord_{n+1}^{fr} \hookrightarrow \Bord_{n+1}^{\partial}.\]
Often the theory $\TFT{Z}^{bulk}$ is required to be invertible and is called an \emph{anomaly}.

Note the dimension: we would like to consider such a theory as a generalization of an ``absolute'' theory of one dimension less valued in the looping $\Omega\mathcal{C}$, namely, the values ``at the free boundary'',
\begin{equation}\label{eqn:deloopTFT}
\TFT{Z}\colon \Bord_n^{fr} \to \Omega\mathcal{C}\,.
\end{equation}

 A slight generalization is a \emph{defect theory}, in which we allow two bulk theories separated by defect, namely a codimension 1 submanifold. In the framed fully extended case, according to the Cobordism Hypothesis with singularities these are characterized by an $n$-adjunctible 1-morphism between two $(n+1)$-dualizable objects in~$\mathcal C$.

In contrast, a \emph{relative or twisted field theory} should not require the bulk theory to be $(n+1)$-dimensional, but rather merely be $n$-dimensional.
Geometrically, this  amounts to only asking for (some) cylinders over the free boundary. Such a cobordism (higher) category was only very recently constructed in William Stewart's PhD thesis \cite{StewartThesis}.

Instead, Stolz--Teichner in \cite{stolz2011supersymmetric} defined a twisted field theory to be a symmetric monoidal natural transformation between symmetric monoidal functors
\begin{equation} \label{eq:TwistedTQFT}
\begin{tikzcd}
   &                                                                     & {} \arrow[dd, "\Longrightarrow" rotate=-90, phantom, shift right] &           \\
\TFT{Z}^{tw}: & \Bord_n \arrow[rr, "S=\uobj", bend left] \arrow[rr, "T"', bend right] &                                       & \ \ca C \,. \ \\
   &                                                                     & {}                                    &          
\end{tikzcd}
\end{equation}
 Here, if $\Bord_n$ is an $(\infty,n)$-category, then $\mathcal{C}$ should be an $(\infty,N)$-category for $N>n$, as above. Moreover, we often require that the source $S=\uobj$ is the trivial functor valued the unit in which case $T$ is
called \emph{twist}.
However, for this definition to be meaningful, the natural transformation must be \emph{lax} or \emph{oplax}, as was explained in the bicategorical situation in \cite{stolz2011supersymmetric}.

For higher categories lax and oplax natural transformations were defined in \cite{JF_CS}. Unfortunately, there is a slightly unsatisfying dichotomy: symmetric monoidal \emph{oplax}  natural transformations from $\uobj$ to $T$ enjoy the desired feature \eqref{eqn:element_in_vsp} that for any closed bordism $b$ we have a morphism $\uobj \to T(b)$, so we are chosing an element in $T(b)$. In the \emph{lax} situation, the arrow switches direction depending on the dimension (see \cite[Example 7.3]{JF_CS}). On the other hand, \emph{lax} symmetric monoidal natural transformations from $\uobj$ to $\uobj$ are equivalent to untwisted (``usual'') field theories valued in $\Omega\mathcal{C}$ and hence satisfy \eqref{eqn:deloopTFT}, whereas \emph{oplax} ones are equivalent to untwisted (``usual'') field theories in $\Omega\mathcal{C}^{\odd \op}$ (\cite[Theorem 7.4 \& Remark 7.5]{JF_CS}).

Furthermore, in \cite{JF_CS}  a classification of framed fully extended twisted field theories was given using the Cobordism Hypothesis in terms of the value at a point. This value is required to be an $a^n$-adjunctible 1-morphism between $n$-dualizable objects, where $a^n=l^n$ or $a^n=r^n$ (see \cref{exam:a-adj}) for lax and oplax, respectively \cite[Theorem 1.6]{JF_CS}. 
We saw in \cref{Coro:evenDimensions} that, as a consequence of our main theorem, in terms of the dualizability/adjunctibility these notions coincide in even dimensions and that in this case there is a second class of conditions we could have chosen. This suggests the question of whether we could have changed the definition of (op)lax TFTs slightly by reversing directions in a way to simultaneously satisfy \eqref{eqn:element_in_vsp} and \eqref{eqn:deloopTFT}.

In this section we prove that none of the choices resolves this dichotomy. Hence, as the first requirement of chosing elements in the theory $T$ seems to be generally considered the more important one, we \emph{define} a relative or twisted field theory to be a symmetric monoidal \emph{oplax} natural transformation from $\uobj$ to some categorified field theory $T$. However, we also prove that any other notion can be obtained by taking oplax natural transformations valued in an appropriate opposite category of the target.

Returning to the geometric picture in \cite{StewartThesis}, these choices correspond to different choices of how the $n$-dimensional framing of the bordism restricts along the defect or boundary. Hence, also geometrically, precisely these choices appear.

Moreover, Freed--Teleman in \cite{FT2014relative} avoid giving a precise definition of relative field theory, and for instance in the situation in \cite{FreedTelemanGapped} a stronger condition than the ones we just considered appears: namely, an $n$-adjunctible 1-morphism between objects which are only $n$-dualizable.

Summarizing, to obtain a framed fully extended boundary/defect, relative/twisted TFT, or dualizable relative field theory, we require the following data for a 1-morphism $f\colon X\to Y$. The above mentioned stronger condition is indicated in the last column.
\smallskip
{\tiny 
\begin{center}
\begin{tabular}{c||c|c|c|c}
&defect  & twisted lax & twisted oplax & \cite{FreedTelemanGapped}\\\hline
objects $X, Y$ &  $(n+1)$-dualizable & $n$-dualizable &  $n$-dualizable &  $n$-dualizable\\
1-morphism $f$  & left \& right adjoints  & left adjoint & right adjoint & left \& right adjoints\\
(co)unit 2-morphisms   & left \& right adjoints & left adjoint & right adjoint & left \& right adjoints\\
	\vdots & \vdots & \vdots & \vdots & \vdots\\
(co)unit $n$-morphisms   & left \& right adjoints & left adjoint & right adjoint & left \& right adjoints\\	
\end{tabular}
\end{center}
\smallskip
}

Our second goal in this section is to compare dualizable relative/twisted field theories to the final option.

\subsection{Comparing lax and oplax natural transformations and variants}

In this section we freely use the notation from \cite{JF_CS}. In brief, given an $(\infty, N)$-category $\mathcal{C}$ two arrow categories are defined, namely $\Cri$ and $\Cdo$, which we will call the \emph{oplax and lax arrow categories of $\mathcal{C}$}, respectively. Both have 1-morphisms in $\mathcal{C}$ as objects and squares as 2-morphisms, but the squares are filled with 2-morphisms in different directions. They are defined by
\[ \CaRight{\ca C}_{\vec k} = \mathrm{map}(\Theta^{(1) ; \vec k}, \ca C)  \qquad \mbox{and} \qquad  \Cdo_{\vec k} = \mathrm{map}(\Theta^{\vec k ; (1)}, \ca C)\]
	and 
	$\Theta^{\vec k ; (1)}$ is glued from certain computads $\Theta^{(j) ; (1)}$, and similarly for $\Theta^{(1) ; \vec k}$, see \cite[Definition 5.7, 5.10]{JF_CS}.

Using these categories, we briefly recall the definition of lax and oplax natural transformations.
\begin{defi}\cite[Definition 1.3 \&  Definition 6.7]{JF_CS}
Let $\mathcal{C}$ be an $(\infty, N)$-category. A \emph{lax natural transformation} is a functor 
\[\mathcal{B} \longrightarrow \Cdo\]
and an \emph{oplax natural transformation} is a functor 
\[\mathcal{B} \longrightarrow \Cri \,.\]
Symmetric monoidal (op)lax natural transformations are symmetric monoidal such functors.
\end{defi}

Using the Cobordism Hypothesis, the following Theorem characterizes (op)lax twisted fully extended framed topological field theories.
\begin{theo}[Theorem 7.6 in \cite{JF_CS}] \label{theo:laxOplaxDualizable}
	Let $n > 0$ and $\ca C$ be a symmetric monoidal 
	$(\infty,N)$-category.
	An object $f:X \to Y$ in $\Cdo$ is $n$-dualizable
	if and only if $X, Y$ are $n$-dualizable and
	$f$ as a morphism in $\ca C$ is $n$-times left adjunctible.
	An object $f:X \to Y$ in $\Cri$ is $n$-dualizable
	if and only if $X, Y$ are $n$-dualizable and
	$f$ as a morphism in $\ca C$ is $n$-times right adjunctible.
\end{theo}

Our first proposition is a connection between the lax and oplax arrow categories.
\begin{prop}\label{thm:Carrows}
	Let $\mathcal{C}$ be an $(\infty, N)$-category. With the opposite function from \cref{ex:odd/even op} we have an equivalence of $(\infty,N-1)$-categories 
	\[\left( \CaRight{{\ca C}^{\odd \op}}\right)^{\odd \op} \simeq \Cdo\]
	which on $j$-morphisms is induced by the isomorphism of computads from \cref{prop_computads}. Moreover, the two maps $s_v,t_v\colon \Cdo \to \ca C$ correspond to the maps
	$$t_h, s_h \colon \left( \CaRight{{\ca C}^{\odd \op}}\right)^{\odd \op} \to {\ca C}^{\odd \op} \,.$$
\end{prop}

\begin{proof}
	We show that
	\[ \CaRight{{\ca C}^{\odd \op}} \simeq \left(\Cdo\right)^{\odd \op}.\]
Since $\Theta^{(1) ; \vec k}$ is glued from $\Theta^{(1) ; (j)}$ and similarly for the indices switched, we obtain isomorphisms
	\[ \left( \Theta^{(1) ; \vec k} \right)^{\odd \op}  \cong \Theta^{\vec k ; (1)} \]
	from the isomorphism of computads from the following \cref{prop_computads} and a tedious matching. 	These isomorphisms induce equivalences of spaces
	\[\left( \CaRight{{\ca C}^{\odd \op}} \right)_{\vec k} \simeq \left(\Cdo\right)_{\vec k} \,. \]
	
	To see that these assemble to an equivalence of higher categories, observe that the face maps in $\CaRight{\ca C}$ and $\Cdo$ are induced by the inclusions
	\[ s_v, t_v \colon \Theta^{(1);(j-1)} \hookrightarrow \Theta^{(1);(j)}  \qquad \mbox{and} \qquad  s_h, t_h \colon \Theta^{(j-1);(1)} \hookrightarrow \Theta^{(j);(1)} ,\]
	respectively.
	In \cref{prop_computads} we identified $s_v, t_v$ on the left with $s_h,t_h$ on the right, the ordering depending on the parity. 
\end{proof}

\begin{lema}\label{prop_computads}
	There is an isomorphism of computads
	\[\left(\Theta^{(1);(j)}\right)^{\odd \op} \cong \Theta^{(j);(1)}, \]
	given as follows. In both cases the top generators exchange, $\theta_{1;j} \leftrightarrow  \theta_{j;1}$. If
	\begin{description}
		\item[$j$ is odd]  
		\begin{align*} 
			s_v\theta_{1;j-1} &\leftrightarrow t_h \theta_{j-1;1} & s_h \theta_{0;j} & \leftrightarrow t_v \theta_{j;0}\\
			t_v\theta_{1;j-1} &\leftrightarrow s_h \theta_{j-1;1} & t_h\theta_{0;j} & \leftrightarrow s_v \theta_{j;0}\,.
		\end{align*}
		\item[$j$ is even]  
		\begin{align*} 
			s_v\theta_{1;j-1} &\leftrightarrow s_h \theta_{j-1;1} & s_h \theta_{0;j} & \leftrightarrow t_v \theta_{j;0}\\
			t_v\theta_{1;j-1} &\leftrightarrow t_h \theta_{j-1;1} & t_h\theta_{0;j} & \leftrightarrow s_v \theta_{j;0}\,.
		\end{align*}
	\end{description}
	Here we use the same notation for the generators of $\Theta_{(1);(j)}^{\odd \op}$ as for $\Theta_{(1);(j)}$.
\end{lema}

\begin{proof}
	The definition of $\theta_{1;j}$ in $\Theta^{(1);(j)}$ says, since 1 is odd,
	\begin{equation}\label{eqn:theta}
		t_v\theta_{1;j-1} \circ s_h \theta_{0;j} \xrightarrow{\theta_{1;j}} t_h \theta_{0;j} \circ s_v\theta_{1;j-1} \,.
	\end{equation}
	If $j$ is odd, since all morphisms appearing in the source and target are $j$-morphisms and hence are reversed in $\left(\Theta^{(1);(j)}\right)^{\odd \op}$, we have in $\left(\Theta^{(1);(j)}\right)^{\odd \op}$ that
	\[s_h \theta_{0;j} \circ t_v\theta_{1;j-1} \xrightarrow{\theta_{1;j}} s_v\theta_{1;j-1} \circ t_h \theta_{0;j}  \,.\]
	
	On the other hand,  if $j$ is odd, in $\Theta^{(j);(1)}$ we have that
	\[t_v \theta_{j;0} \circ s_h\theta_{j-1;1} \xrightarrow{\theta_{j;1}} t_h\theta_{j-1;1} \circ s_v \theta_{j;0}  \,.\]
	Comparing the last two expressions we find the exchange.
	
	If $j$ is even, in \eqref{eqn:theta} the morphisms appearing in the source and target are $j$-morphisms and hence are not reversed in $\left(\Theta^{(1);(j)}\right)^{\odd \op}$, but $\theta_{1;j}$ is a $(j+1)$-morphism and hence reversed, so the source and target switch:
	\[t_h \theta_{0;j} \circ s_v\theta_{1;j-1}  \xrightarrow{\theta_{1;j}}t_v\theta_{1;j-1} \circ s_h \theta_{0;j}\,.\]
	
	On the other hand, in $\Theta^{(j);(1)}$ we have that
	\[s_v \theta_{j;0} \circ s_h\theta_{j-1;1} \xrightarrow{\theta_{j;1}} t_h\theta_{j-1;1} \circ t_v \theta_{j;0}  \,.\]
	
	Now proceed by induction.
\end{proof}

Inspired by \cref{thm:Carrows}, we generalize the (op)lax arrow categories to general opposite functions. Recall from \cref{Prop:oppositeFunctionAdjunction} the opposite function $\op_{a^n,b^n,k}^N$ associated to a pair of dexterity functions. In this section, we will always use $b^n=r^n$ and $k=1$, hence we abbreviate notation to
\[
	\op_{a^n} \coloneqq \op_{a^n,r^n,1}.
\]

\begin{defi}
	Let $\ca C$ be a $\infN N$category and $a^n$ a dexterity function.
	Define the \emph{$a^n$-lax arrow category} to be the $(\infty, N-1)$-category
	\[
		{\ca C}^{a^n} := \left( \CaRight{ {\ca C}^{\op_{a^n}} }\right)^{\op_{a^n}} \,.
	\]
\end{defi}

We could also have used the other opposite function from \cref{Prop:oppositeFunctionAdjunction} in the definition above
and we would get an equivalent arrow category.
To see this we use the following Proposition that we obtain similarly to \cref{thm:Carrows}.
\begin{prop} \label{Prop:opppCommutesWithRightArrow} 
	Let $\mathcal{C}$ be an $(\infty, N)$-category. Let $\mathrm{oppp}  \equiv \op$ be the constant opposite function valued $\op$. Then there is an equivalence of $(\infty,N-1)$-categories
	\[
		\left( \CaRight{ {\ca C}^{\oppp} }\right)^{\oppp} \simeq \Cri \,
	\]
	induced on $j$-morphisms by  isomorphisms of computads as indicated in \eqref{eqn:theta_oppp}.
	
	Moreover, the two maps  $s_h, t_h \colon\Cri \to \ca C$ correspond to the maps
	$$t_h, s_h \colon \left( \CaRight{ {\ca C}^{\oppp} }\right)^{\oppp} \to {\ca C}^{\oppp} \,.$$
	
\end{prop}
\begin{proof}
	Recall that in $\Theta^{(1);(j)}$, we have that
	\begin{equation*}
		t_v\theta_{1;j-1} \circ s_h \theta_{0;j} \xrightarrow{\theta_{1;j}} t_h \theta_{0;j} \circ s_v\theta_{1;j-1} \,.
	\end{equation*}
	Hence, in $(\Theta^{(1);(j)})^\oppp$ we have that
	\[
		s_h \theta_{0;j}\circ t_v\theta_{1;j-1} \xleftarrow{\theta_{1;j}} s_v\theta_{1;j-1} \circ t_h \theta_{0;j} \,.
	\]
	Switching $s_v$ and $t_v$, and replacing $s_h$ and $t_h$, we obtain an isomorphism of computads
	\begin{equation}\label{eqn:theta_oppp}
		(\Theta^{(1);(j)})^\oppp \cong (\Theta^{(1);(j)})^\oppp \,.
	\end{equation}
\end{proof}

Since $\mathcal{D}^{-\op_{a^n,r^n,1}} = {\left(\mathcal{D}^{\op_{a^n,r^n,1}}\right)}^{\oppp}$ holds for any $\infN M$category $\ca D$
we immediately conclude the following.
\begin{coro} \label{Coro:OppositeFunctionChoiceDoesntMatter}
	Let $\ca C$ be an $\infN N$category and $a^n$ a dexterity function.
	Then there is an equivalence of $(\infty,N-1)$-categories
	\[
		{\ca C}^{a^n} \simeq  \left( \CaRight{ {\ca C}^{-\op_{a^n}} }\right)^{-\op_{a^n}}.
	\]
\end{coro}

\begin{exam}
	If $a^n=r^n$ we have from \cref{exam:OddOppositeFunction} that $\op_{a^n} \equiv \id$ is the constant opposite function valued $\id$.
	In this case we have ${\ca C}^{a^n} = \CaRight{\ca C}$.
	\exampleend
\end{exam}

\begin{exam}
	If $a^n=l^n$ and $n+1=N$ we seen in \cref{exam:OddOppositeFunction} that $\op_{a^n}=\even \op$.
	Hence ${\ca C}^{a^n} = { \left( \CaRight{ {\ca C}^{\even \op} } \right) }^{\even \op}$ and by \cref{Coro:OppositeFunctionChoiceDoesntMatter}
	and \cref{thm:Carrows} we have
	 \[
	 	{\ca C}^{a^n} \simeq { \left( \CaRight{ {\ca C}^{\odd \op} } \right) }^{\odd \op}
	 	\simeq \Cdo.
	 \]
	 \exampleend
\end{exam}

Theorem \ref{theo:laxOplaxDualizable} generalizes to these generalized arrow categories.
\begin{coro}\label{cor:a-lax_dualizability}
	Let $\ca C$ be a symmetric monoidal $\infN N$category, $f:X \to Y$ a $1$-morphism in $\ca C$, and $a^n$ a dexterity function. The following  are equivalent:
	\begin{enumerate}
	\item 
	When viewed as an object in ${\ca C}^{a^n}$, $f$ is $n$-dualizable.
	\item The objects $X$ and $Y$ are $n$-dualizable and $f$ is $n$-times right adjunctible in~${\ca C}^{\op_{a^n}}$.
	\item The objects $X$ and $Y$ are $n$-dualizable and $f$ is $a^n$-adjunctible in~$\ca C$.
	\end{enumerate}
\end{coro}
\begin{proof}
	Let $\op_{a^n} = \op_{a^n,r^n,1}$ be the opposite function associated to $a^n$ and $r^n$ from \cref{Prop:oppositeFunctionAdjunction}.
	Since taking opposites does not change dualizability, the 1-morphism $f$ is dualizable in ${\ca C}^{a^n}$ if and only if it is dualizable in $\CaRight{{\ca C}^{\op_{a^n}}}$.
	This is the case if and only if $X$ and $Y$ are $n$-dualizable in ${\ca C}^{\op_{a^n}}$ 
	and $f$ is $n$-times right adjunctible in ${\ca C}^{\op_{a^n}}$, 
	by applying \cref{theo:laxOplaxDualizable} to $\CaRight{{\ca C}^{\op_{a^n}}}$.
	By \cref{Prop:oppositeFunctionAdjunction} the latter is equivalent to $f$ being $a^n$-adjunctible
	in $\ca C$.
\end{proof}

Our conclusion now is the following: we \emph{could} define an $a^n$-lax twisted field theory to be a symmetric monoidal ``$a^n$-lax natural transformation'', i.e.~a symmetric monoidal functor
\[\Bord_n \longrightarrow {\ca C}^{a^n},\]
and \cref{cor:a-lax_dualizability} characterizes the framed fully extended ones in terms of $a^n$-adjunctability.
With this at hand, we have a (many!) notion(s) of $a^n$-lax twisted field theory which in terms of dualizability lies in the second equivalence of dexterity functions. However, this does not resolve the dichotomoy between \eqref{eqn:element_in_vsp} and \eqref{eqn:deloopTFT}:

If we want to satisfy \eqref{eqn:element_in_vsp}, we must choose the \emph{oplax} case, since the other choices involve taking opposites at certain levels. No other dexterity function can achieve \eqref{eqn:element_in_vsp}.

As for the other desideratum \eqref{eqn:deloopTFT}, a variant of \cite[Theorem 7.4 \& Remark 7.5]{JF_CS} shows that since we are taking opposites at various levels, this is only satisfied for the \emph{lax} case.

\subsection{Dualizable relative topological field theories}
We now turn to comparing the adjunctibility conditions in the table in the beginning of this section. Leaving the first column aside, observe that by ambidexterity, i.e.~\cref{Theo:fromMixedToNadjunctibility} (a) the difference only lies in the adjunctibility of the appearing (co)unit $n$-morphisms.
Hence, we only need to check that the appearing (co)unit $n$-morphisms also have the other adjoint.
This can be guaranteed by asking for the twisted/relative field theory to be itself adjunctible, as was suggested to us by Constantin Teleman.

In the following propositions we assume that the Cobordism Hypothesis holds, that $\ca C$ is an $\infN N$category and $n \geq 1$.
\begin{prop}\label{prop:dualizable_relative}
Let $\TFT{Z}^{tw}$ be a framed fully extended $n$-dimensional oplax twisted field theory, i.e.~a symmetric monoidal oplax natural transformation as above. 	If $\TFT{Z}^{tw}$, when viewed as a 1-morphism in the $(\infty,n)$-category $\mathrm{Fun}^{oplax}(\Bord_n^{fr}, \mathcal{C})$, is adjunctible, then the value of $\TFT{Z}^{tw}$ at a point is a 1-morphism $f\colon X\to Y$ in $\mathcal{C}$ such that $X$ and $Y$ are $n$-dualizable and $f$ is $n$-adjunctible. 

The same statement is true for ``oplax'' replaced by ``lax''.
\end{prop}

\begin{proof}
By the characterization \cite[Theorem 1.6]{JF_CS} the value of $\TFT{Z}^{tw}$ at a point is a 1-morphism $f\colon X\to Y$ in $\mathcal{C}$ such that $X$ and $Y$ are $n$-dualizable and $f$ is $n$-times right adjunctible. 

Evaluating the left adjoint $\TFT{Z}^L$ of $\TFT{Z}^{tw}$, we obtain a 1-morphism $g\colon Y\to X$ in $\mathcal{C}$, which again by \cite[Theorem 1.6]{JF_CS}  is $n$-times right adjunctible.
Moreover, evaluating the adjunction data of the adjunction $\TFT{Z}^{L} \adj \TFT{tw}^L$ at a point, we see that $g$ is a left adjoint of $f$. By \cref{Coro:adjointsOfNadjunctible} $f$ is n-adjunctible.	
\end{proof}

\begin{prop} 
	Let  $j \geq 1$. 
	Let $\TFT{Z}^{otw}$ be a framed fully extended $n$-dimensional oplax twisted field theory and the 1-morphism $f\colon X \to Y$ in $\ca C$ its value at a point.
	Then, the following are equivalent:
	\begin{enumerate}
		\item In $\mathrm{Fun}^{oplax}(\Bord_n^{fr}, \mathcal{C})$ we have that the 1-morphism $\TFT{Z}^{otw}$ is $\even^j$-adjunctible.
		\item In $\mathcal{C}$ we have that $X$ and $Y$ are $n$-dualizable and $f$ is $\even^{n+j}$-adjunctible.
	\end{enumerate}
	
	Let $\TFT{Z}^{ltw}$ be a framed fully extended $n$-dimensional lax twisted field theory and the 1-morphism $g\colon X \to Y$ in $\ca C$ its value at a point.
	Then, the following are equivalent:
	\begin{enumerate}
		\item In $\mathrm{Fun}^{lax}(\Bord_n^{fr}, \mathcal{C})$ we have that the 1-morphism $\TFT{Z}^{ltw}$ is $\even^j$-adjunctible.
		\item In $\mathcal{C}$ we have that $X$ and $Y$ are $n$-dualizable and 
		\begin{align*}		 
			&g \ \text{is} \
			\begin{cases}
				\even^{n+j}\text{-adjunctible} & \text{if } n \text{ is even}\,, \\
				\odd^{n+j}\text{-adjunctible} & \text{if } n \text{ is odd}\,.
			\end{cases}			
		\end{align*}
	\end{enumerate}
The same is true with all appearances of ``even'' and ``odd'' reversed.
\end{prop}

\begin{proof}
	We consider the oplax case first.
	$(1)\Rightarrow (2)$
	If $\TFT{Z}^{otw}$ is $\even^j$-adjunctible, it has a right adjoint $\TFT{Z}^R$. Evaluating $\TFT{Z}^R$ at a point we obtain a 1-morphism $f^R$ that is a right adjoint of $f$. The 
	unit and counit of the adjunction $\TFT{Z}^{otw} \adj \TFT{Z}^R$ are 2-morphisms $\TFT{U}$ and $\TFT{C}$ in $\mathrm{Fun}^{oplax}(\Bord_n^{fr}, \mathcal{C})$. By \cite[Corollary 7.7]{JF_CS} evaluating $\TFT{U}$ and $\TFT{C}$ at a point we obtain 2-morphisms $u$ and $c$ in $\ca C$ which are $n$-times right adjunctible.
	Hence, $f$ is $(n+1)$-times right adjunctible.
	
	By induction, by iterating this argument for the adjunctibility data of the adjunction $\TFT{Z}^{otw} \adj \TFT{Z}^R$,  we have that this evaluates at a point to adjunctibility data exhibiting $f$ as $(n+j)$-times right adjunctible, which is the same as $\even^{n+j}$-adjunctible.
	
	If $\TFT{Z}^{otw}$ is $\odd^j$-adjunctible we can use the same induction to get 
	$(n+j-1)$-times right adjunctibility data for $f$. However, in the last step instead of right adjoints for the unit and counit $j$-morphisms of the data of $\TFT{Z}^{otw}$, we have left adjoints. After evaluation at a point we obtain
	left adjoints of the unit and counit $j$-morphisms of the $(j-1)$-times right adjunctibility data for $f$ that are themselves $n$-times right adjunctible by
	\cite[Corollary 7.7]{JF_CS}. Combining these we get
	$a^{n+j}$-adjunctibility, with
	\begin{equation} \label{Eq:oddDexterityWithJmapsToL}
		a^{n+j}: \{ 1, \ldots, n + j \} \to \{ L, R \}, \quad i \mapsto 
		\begin{cases}
			L & i = j \,,  \\
			R & i \neq j
		\end{cases}
	\end{equation}
	which by \cref{Prop:AnReducesToTwoClasses} is equivalent to $\odd^{n+j}$-adjunctibility data for $f$.
	
	$(2)\Rightarrow (1)$ 
	Let $f\colon X \to Y$ be an $\even^{n+j}$-adjunctible. 1-morphism between $n$-dualizable objects. Then $f$ is $n$-adjunctible by \cref{Theo:fromMixedToNadjunctibility}. Moreover, $f$ has a right adjoint $f^R$ which itself is $n$-adjunctible by \cref{Prop:AdjointAlsoAdjunctible}, and in particular $n$-times right adjunctible. Hence, by \cite[Theorem 1.6]{JF_CS} we obtain a framed fully extended $n$-dimensional oplax twisted field theory $\TFT{Z}^R$. We claim that this is a right adjoint.
	
	The unit $u$ and counit $c$ of $f\adj f^R$ are $n$-times right adjunctible because $f$ is $(n+1)$-times right adjunctible. So,  invoking \cite[Corollary 7.7]{JF_CS} there are 2-morphisms $\TFT{U}$ and $\TFT{C}$ in $\mathrm{Fun}^{oplax}(\Bord_n^{fr}, \mathcal{C})$ whose value at a point are $u$ and $c$, respectively. We claim that they are the unit and counit of an adjunction. To see this, we check the compositions in the \eqref{eq:zig} and \eqref{eq:zag} at the point, where they are identities because $u$ and $c$ are the (co)units of an adjunction. Hence, by \cite[Corollary 7.7]{JF_CS} they are identities in $\mathrm{Fun}^{oplax}(\Bord_n^{fr}, \mathcal{C})$.
	
	We now repeat the argument inductively, replacing $f\adj f^R$ by an adjunction appearing in the adjunctibility data of $f$ to produce adjunctibility data of $\TFT{Z}^{otw} \adj \TFT{Z}^R$ exhibiting $j$-times right adjunctibility.
	
	If $f\colon X \to Y$ is $\odd^{n+j}$-adjunctible we can
	construct $(j-1)$-times right adjunctibility data for $\TFT{Z}^{otw}$
	as above. For the last step we use that $f$ by \cref{Prop:AnReducesToTwoClasses} is also $a^{n+j}$-adjunctible as defined in \eqref{Eq:oddDexterityWithJmapsToL}.
	From the $a^{n+j}$-adjunctibility data we get left adjoints of the unit and counit $j$-morphisms that are themselves $n$-times right adjunctible. By \cite[Corollary 7.7]{JF_CS} and arguing as above we complete the $(j-1)$-times right adjunctibility data for $\TFT{Z}^{otw}$ to $\odd^{j}$-adjunctibility data.

	The lax case works similar with appropriate dexterity functions given by \cref{Prop:AnReducesToTwoClasses}.
\end{proof}

Combining the Proposition above with \cref{Theo:fromMixedToNadjunctibility} we obtain
the following.
\begin{coro}\label{coro:characterization_dualizable_twisted}
	Let  $j \geq 1$. Let $\TFT{Z}^{tw}$ be a framed fully extended $n$-dimensional oplax twisted field theory and the 1-morphism $f\colon X \to Y$ in $\ca C$ its value at a point. Then, the following are equivalent:
	\begin{enumerate}
		\item In $\mathrm{Fun}^{oplax}(\Bord_n^{fr}, \mathcal{C})$ we have that the 1-morphism $\TFT{Z}^{tw}$ is $j$-adjunctible.
		\item In $\mathcal{C}$ we have that $X$ and $Y$ are $n$-dualizable and $f$ is $(n+j)$-adjunctible. 
	\end{enumerate}
	The same statement is true for ``oplax'' replaced by ``lax''.
\end{coro}

\section{More general notions of higher adjunctibility and binary trees} \label{sec:moreGeneral}

In this section we go one step further in generalizing the definition of mixed higher adjunctibility.
The following example illustrates that in a higher Morita category we can have a kind of $2$-adjunctibility that fits neither of the $4$ possible dexterity functions `$LL$', `$RR$'
`$RL$' or `$LR$'.

\begin{exam} \label{Ex:moritaTree}
	Let $k$ be a field or a commutative ring. 
	We have seen a description of the Morita $3$-category $\Alg 2$ in \cref{sec:sec2}.
	
	Let $_k A_k$ be a $1$-morphism $k \to k$ in $\Alg 2$.
	As discussed in \cref{sec:sec2} the right and left adjoints are given by the opposite algebra
	${}^{\vphantom{\op}}_k A^{\op}_k$.
	For the adjunction $_k A_k \adj {}^{\vphantom{\op}}_k A^{\op}_k$ we have 
	\begin{align*}
		& \phantom{\text{co}}\text{unit} \qquad 
		\, u : \id_k \to { _k A_k } \, {\te_k} \, {}^{\vphantom{\op}}_k A^{\op}_k \,  \ 
		= \  _{ _k k_k } A_{ ( { _k A_k } ) { \, {\te_k} \, } ( { {}^{\vphantom{\op}}_k A^{\op}_k } ) } 
		\qquad \text{and}
		\\ \\ & \text{counit}
		\qquad  \, c : { {}^{\vphantom{\op}}_k A^{\op}_k } \, {\te_k} \, _k A_k \to \, \id_k \,  \ 
		= \  _{ ( { {}^{\vphantom{\op}}_k A^{\op}_k } ) \, {\te_k} \, ( { _k A_k } ) } A_{ _k k_k }
		\, .
	\end{align*}
	
	By \cref{Prop:bimodDual} the unit $u$ is left adjunctible if and only if $A$ is a finite dimensional $k$-vector space and $u$ is right adjunctible if and only if $A$ is finitely presented projective as a $A \te_k A^{\op}$ module. Since a finitely presented module is projective if and only if it is flat, the latter condition is equivalent to $A$ being a separable $k$-algebra.
	Similarly, the counit $c$ is right adjunctible if and only if $A$ is finite dimensional and 
	it is left adjunctible if and only if $A$ is separable.
	
	Choosing the dual numbers $A = k[x]/(x^2)$ for $A$ we have an example where $A$ is finite dimensional over $k$ but not separable.
	In this case the $1$-morphism $_k A_k$ does not satisfy any $a^2$-adjunctibility since the unit is left but not right adjunctible, while the counit is right but not left adjunctible.
	\exampleend
\end{exam}

To capture this phenomenon we introduce dexterity trees, that
specify for each unit and counit individually whether
it should have a left or right adjoint.
The tree structure arises from the form of our data. For example, in \ref{Ex:moritaTree} with $A = k[x]/(x^2)$ we have adjunctions
\begin{displaymath}
	\vcenter{\hbox{\begin{tikzpicture}
				\tikzstyle{level 1}=[level distance=8mm,sibling distance=30mm]
				\node {$(f^L \adj f, u, c)$}
				child {node {$(u^L \adj u, u_u, c_u)$}
				}
				child {node {$(c \adj c^R, u_c, c_c)$}
				};
	\end{tikzpicture}}}
	\quad
	\text{with dexterity tree}
	\quad
	\vcenter{\hbox{\begin{tikzpicture}
				\tikzstyle{level 1}=[level distance=8mm,sibling distance=10mm]
				\node {$L$}
				child {node {$L$}
				}
				child {node {$R$}
				};	
	\end{tikzpicture}}}
\end{displaymath}
Before we get to the precise Definition in \ref{Defi:treeAdj} we recall some vocabulary regarding trees.

The complete rooted binary tree $\tree {n+1}$ of height $n+1$ is given by words on the alphabet $\{0,1\}$ of length $\leq n$. The \emph{root} and is  $\emptyset$ and the \emph{parent} of a word $w$ of length $l$ is the length $l-1$ prefix of $w$.
We refer to a direct successor of a node as \emph{child} and to two children with the same parent as \emph{siblings}.
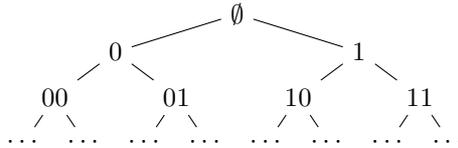
\begin{center}
	\begin{tikzpicture}
		\tikzstyle{level 1}=[level distance=5mm,sibling distance=32mm]
		\tikzstyle{level 2}=[level distance=6mm,sibling distance=16mm]	
		\tikzstyle{level 3}=[level distance=6mm,sibling distance=8mm]				
		\node{$\emptyset$}
		child {node {$0$}
			child {node {$00$}
				child {node {$\cdots$}}
				child {node {$\cdots$}}
			}
			child {node {$01$}
				child {node {$\cdots$}}
				child {node {$\cdots$}}
			}
		}
		child {node {$1$}
			child {node {$10$}
				child {node {$\cdots$}}
				child {node {$\cdots$}}
			}
			child {node {$11$}
				child {node {$\cdots$}}
				child {node {$\cdots$}}
			}
		};
	\end{tikzpicture}
\end{center}
The length $|w|$ of a word $w$ is the \emph{depth} of the corresponding node in the tree.
Note that we use the convention that the height of a rooted tree is the number of nodes in the longest downward path starting at the root. Thus $\tree n$ has $n$ layers and $2^{n}-1$ nodes.
This makes the number of entries $n$ of a dexterity function $a^n$ align with the height of the corresponding dexterity tree from \cref{Defi:treeAdj}.
In this section, by tree we always mean a complete rooted binary tree of finite height.

\begin{defi} \label{Defi:treeAdj}
	Let $\ca C$ be an $\infN N$category.
	Let 
	\[
	t^n: \tree n \to \{ L, R \}.
	\]
	be a complete rooted binary tree labelled by $\{ L, R \}$ of height $n$.
	Let $f$ be a $k$-morphism in $\ca C$.
	A set of \emph{$t^n$-adjunctibility data} for $f$ with \emph{dexterity tree $t^n$} is defined inductively as
	\begin{itemize}[itemsep=2mm]
		\item step $1$: the data of a
		$\begin{cases}
			\text{left adjoint, } & \text{if} \quad t^n(\emptyset) = L, \\
			\text{right adjoint, } & \text{if} \quad t^n(\emptyset) = R,
		\end{cases}$ \\
		together with unit and counit witnessing the adjunction.
		\item step $j= 2, \ldots, n$: for
		all units from step $j-1$, witnessing the adjunction corresponding to the length $j-1$ word $w$,
		\\ \phantom{step $1$:} the data of a
		$\begin{cases}
			\text{left adjoint, } & \text{if} \quad t^n(w \, 0) = L, \\
			\text{right adjoint, } & \text{if} \quad t^n(w \, 0) = R,
		\end{cases}$ \\
		together with units and counits witnessing the adjunctions.
		Similarly, data for all counits depending on the value of $t^n(w \, 1)$.
	\end{itemize}
	A $k$-morphism $f$ is \emph{$t^n$-adjunctible} if there exists a set of 
	$t^n$-adjunctibility data for~$f$.
\end{defi}

Corresponding to \cref{Lemma:GeneralInterchangeLemma},
that allowed us to exchange `$LL$' $\sim$ `$RR$' and
`$LR$' $\sim$ `$RL$' in the entries of dexterity functions to get equivalent notions of $a^n$-adjunctibility,
we now can exchange any subtrees of the form
\begin{equation} \label{Eq:subtreeChange}
	\vcenter{\hbox{\begin{tikzpicture}
				\tikzstyle{level 1}=[level distance=6mm]
				\tikzstyle{level 2}=[level distance=7mm,sibling distance=12mm]
				\tikzstyle{level 3}=[level distance=7mm,sibling distance=5mm]				
				\node {}
				child {node {$L$}
					child {node {$L$}
						child {node {$t_1$}}
						child {node {$t_2$}}
					}
					child {node {$L$}
						child {node {$t_3$}}
						child {node {$t_4$}}
					}
				};
	\end{tikzpicture}}}
	\sim
	\vcenter{\hbox{\begin{tikzpicture}
				\tikzstyle{level 1}=[level distance=6mm]
				\tikzstyle{level 2}=[level distance=7mm,sibling distance=12mm]
				\tikzstyle{level 3}=[level distance=7mm,sibling distance=5mm]				
				\node {}
				child {node {$R$}
					child {node {$R$}
						child {node {$t_3$}}
						child {node {$t_4$}}
					}
					child {node {$R$}
						child {node {$t_1$}}
						child {node {$t_2$}}
					}
				};
	\end{tikzpicture}}}
	\quad \text{and} \quad
	\vcenter{\hbox{\begin{tikzpicture}
				\tikzstyle{level 1}=[level distance=6mm]
				\tikzstyle{level 2}=[level distance=7mm,sibling distance=12mm]
				\tikzstyle{level 3}=[level distance=7mm,sibling distance=5mm]				
				\node {}
				child {node {$L$}
					child {node {$R$}
						child {node {$t_1$}}
						child {node {$t_2$}}
					}
					child {node {$R$}
						child {node {$t_3$}}
						child {node {$t_4$}}
					}
				};
	\end{tikzpicture}}}
	\sim
	\vcenter{\hbox{\begin{tikzpicture}
				\tikzstyle{level 1}=[level distance=6mm]
				\tikzstyle{level 2}=[level distance=7mm,sibling distance=12mm]
				\tikzstyle{level 3}=[level distance=7mm,sibling distance=5mm]				
				\node {}
				child {node {$R$}
					child {node {$L$}
						child {node {$t_3$}}
						child {node {$t_4$}}
					}
					child {node {$L$}
						child {node {$t_1$}}
						child {node {$t_2$}}
					}
				};
	\end{tikzpicture}}}
	.
\end{equation}
to get equivalent notions of $t^n$-adjunctibility.
Note that we swapped the subtrees of the two children since \cref{Lemma:MainLemma}
exchanges the role of the unit and counit adjunctions in the data. 
The relations \eqref{Eq:subtreeChange} generate an equivalence relation $\sim$ on the set of dexterity trees.

\begin{theo} \label{Theo:TreeAdjunctibilityEquivalence}
	Let $\ca C$ be an $\infN N$category and $t^n, s^n: \tree n \rightrightarrows \{ L, R \}$ be
	two dexterity trees that are equivalent with respect to $\sim$. Then a $k$-morphism $f$ is {$t^n$-adjunctible} if and only if it is {$s^n$-adjunctible}.
\end{theo}
\begin{proof}
	If $t^n \sim s^n$ are equivalent, there is a series of exchanges of the form \eqref{Eq:subtreeChange} to
	transform $t^n$ into $s^n$. Given a $k$-morphism $f$, all of the intermediate dexterity trees as well as $t^n$ and $s^n$ require
	equivalent adjunctibility data for $f$ by \cref{Lemma:MainLemma}.
\end{proof}

\begin{exam}
	For $3$-times left adjunctibility given by the constant `$L$' dexterity function $l^3$ we have the corresponding constant `$L$' dexterity tree of height $3$. Some equivalent trees are
	\begin{displaymath}
		\vcenter{\hbox{\begin{tikzpicture}
					\tikzstyle{level 1}=[level distance=6mm,sibling distance=12mm]
					\tikzstyle{level 2}=[level distance=6mm,sibling distance=6mm]			
					\node {$L$}
					child {node {$L$}
						child {node {$L$}}
						child {node {$L$}}
					}
					child {node {$R$}
						child {node {$R$}}
						child {node {$R$}}
					};
		\end{tikzpicture}}}
		\mkern9mu , \quad
		\vcenter{\hbox{\begin{tikzpicture}
					\tikzstyle{level 1}=[level distance=6mm,sibling distance=12mm]
					\tikzstyle{level 2}=[level distance=6mm,sibling distance=6mm]			
					\node {$L$}
					child {node {$R$}
						child {node {$R$}}
						child {node {$R$}}
					}
					child {node {$R$}
						child {node {$R$}}
						child {node {$R$}}
					};
		\end{tikzpicture}}}
		\mkern9mu , \quad
		\vcenter{\hbox{\begin{tikzpicture}
					\tikzstyle{level 1}=[level distance=6mm,sibling distance=12mm]
					\tikzstyle{level 2}=[level distance=6mm,sibling distance=6mm]			
					\node {$R$}
					child {node {$R$}
						child {node {$L$}}
						child {node {$L$}}
					}
					child {node {$R$}
						child {node {$L$}}
						child {node {$L$}}
					};
		\end{tikzpicture}}}
		\quad \text{and} \quad
		\vcenter{\hbox{\begin{tikzpicture}
					\tikzstyle{level 1}=[level distance=6mm,sibling distance=12mm]
					\tikzstyle{level 2}=[level distance=6mm,sibling distance=6mm]			
					\node {$R$}
					child {node {$L$}
						child {node {$R$}}
						child {node {$R$}}
					}
					child {node {$R$}
						child {node {$L$}}
						child {node {$L$}}
					};
		\end{tikzpicture}}}
		.
	\end{displaymath}
	The second and third tree are induced by the dexterity functions `$LRR$' and `$RRL$', while the first and fourth tree are not induced by any dexterity function.
	\exampleend
\end{exam}

We saw that for dexterity functions and $a^n$-adjunctibility the a priori $2^n$ different
definitions reduced to two equivalence classes represented by $\even^n$ and $\odd^n$ defined in \eqref{eq:evenN} and \eqref{eq:oddN}.
For dexterity trees and $t^n$-adjunctibility we start with $2^{2^n-1}$ a priori
different notions.
Let $\treeSet{n}$ be the set of dexterity trees and $\treeEqSet{n} := \treeSet{n} / \sim$ be the set of equivalence classes 
generated by the relations \eqref{Eq:subtreeChange}.
The question of determining an upper bound for non-equivalent definitions of dexterity trees and $t^n$-adjunctibility amounts to a
calculation of $|\treeEqSet{n}|$, which is one of the goals for the remainder of this section.
The second one is to find representatives that replace $\even^n$ and $\odd^n$ in this more general setup.
We start with the latter and we will achieve the former goal by counting the representatives.

First we give a name to those $t^n$, that have no other equivalent dexterity trees.
\begin{defi}
	Let $t^n$ be a dexterity tree. We say $t^n$ is a \emph{fixed tree} if every pair of siblings in $t^n$ has distinct labels.
\end{defi}
In a fixed tree we cannot apply any of the relations \eqref{Eq:subtreeChange} to any subtree.

\begin{defi}
	Let $t^n$ be a dexterity tree. We say $t^n$ is in \emph{normal form} if every subtree of $t^n$ is either a fixed tree or has a root $r$ with label `$R$'.
\end{defi}
Any tree in normal form is a constant `$R$' tree
with some subtrees replaced by fixed trees.

\begin{theo} \label{Theo:TreeTheorem}
	Every dexterity tree $t^n$ is equivalent to exactly one tree in normal form $\tau^n$.
	In particular, we have
	\begin{displaymath}
		|\treeEqSet{1}| = 2 
		\quad \text{and} \quad
		|\treeEqSet{n}| = |\treeEqSet{n-1}|^2 + 2^{2^{n-1} -1 }
		\quad \text{for} \quad n \geq 2.
	\end{displaymath}
\end{theo}
We will prove this theorem in a series of Lemmas. 
But before that let us take a look at the example of trees of height $2$.

\begin{exam}
	For height $n=2$ we have $8$ different dexterity trees that contain two equivalent pairs. 
	
	\begin{center}
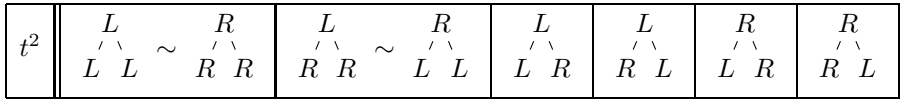
		\begin{tabular}{ | c || l | r | l | r | c | c | c | c |  } \hline
			\raisebox{4mm}{$t^2$}
			& 
			\multicolumn{2}{c|}{		
			\begin{tikzpicture}
				\tikzstyle{level 1}=[level distance=6mm,sibling distance=5mm]
				\node {$L$}
				child {node {$L$}}
				child {node {$L$}};	
			\end{tikzpicture}
			\raisebox{4mm}{$\sim$}
			\begin{tikzpicture}
				\tikzstyle{level 1}=[level distance=6mm,sibling distance=5mm]
				\node {$R$}
				child {node {$R$}}
				child {node {$R$}};	
			\end{tikzpicture}
			} 
			& 
			\multicolumn{2}{c|}{
			\begin{tikzpicture}
				\tikzstyle{level 1}=[level distance=6mm,sibling distance=5mm]
				\node {$L$}
				child {node {$R$}}
				child {node {$R$}};	
			\end{tikzpicture}	
			\raisebox{4mm}{$\sim$}
			\begin{tikzpicture}
				\tikzstyle{level 1}=[level distance=6mm,sibling distance=5mm]
				\node {$R$}
				child {node {$L$}}
				child {node {$L$}};	
			\end{tikzpicture}
			}	
			&
			\begin{tikzpicture}
				\tikzstyle{level 1}=[level distance=6mm,sibling distance=5mm]
				\node {$L$}
				child {node {$L$}}
				child {node {$R$}};	
			\end{tikzpicture}
			& 
			\begin{tikzpicture}
				\tikzstyle{level 1}=[level distance=6mm,sibling distance=5mm]
				\node {$L$}
				child {node {$R$}}
				child {node {$L$}};	
			\end{tikzpicture}
			& 
			\begin{tikzpicture}
				\tikzstyle{level 1}=[level distance=6mm,sibling distance=5mm]
				\node {$R$}
				child {node {$L$}}
				child {node {$R$}};	
			\end{tikzpicture}
			& 
			\begin{tikzpicture}
				\tikzstyle{level 1}=[level distance=6mm,sibling distance=5mm]
				\node {$R$}
				child {node {$R$}}
				child {node {$L$}};	
			\end{tikzpicture}
			\\ \hline
		\end{tabular}
	\end{center}
	Only the first and third dexterity tree are not in normal form, but they are equivalent to one in normal form.
	\exampleend
\end{exam}

\begin{lema} \label{Lemma:canFindNormalForm}
	Let $t^n$ be a dexterity tree. Then there is a dexterity tree $\tau^n$ in normal form with $t^n \sim \tau^n$.
\end{lema}
\begin{proof}
	We can transform a dexterity tree $t^n$ into its equivalent normal form by calling the function
	\Call{TransferIntoNormalform}{$t^n$} defined in algorithm \ref{Alg:TransformIntoNormalform}.
\end{proof}

\begin{algorithm}[ht]
	\caption{Change the label of the root in a tree $t$ if $t$ is not a fixed tree}
	\label{Alg:ChangeRoot}
	\begin{algorithmic}[1] 
		\Function{ChangeRoot}{tree $t$}
		\If{$t$ is a fixed tree}
		\State \Return $t$
		\Else 
		\State Find a node $N$ of minimal depth, which has a sibling with the same label	
		\While{$N$ is not the root of $t$}
		\State Apply suitable \eqref{Eq:subtreeChange} to the parent of $N$ and its children in $t$
		\State $N$ $\gets$ parent of $N$
		\EndWhile
		\State \Return $t$
		\EndIf
		\EndFunction
	\end{algorithmic}
\end{algorithm}

\begin{algorithm}[ht]
	\caption{Change a tree $t$ into normal form}
	\label{Alg:TransformIntoNormalform}
	\begin{algorithmic}[1] 
		\Function{TransformIntoNormalform}{tree $t$}
		\If{$t$ is a fixed tree}
		\State \Return $t$
		\Else 
		\If{the root of $t$ has label `$L$'}
		\State $t$ $\gets$ \Call{ChangeRoot}{$t$} \Comment{(see algorithm \ref{Alg:ChangeRoot})}
		\EndIf
		\State tree $t_{l \phantom{r}}$ := left subtree of the root of $t$
		\State tree $t_{r \phantom{l}}$ := right subtree  of the root of $t$
		\State $t_{l \phantom{r}}$ $\gets$ \Call{TransferIntoNormalform}{$t_l$}
		\State $t_{r \phantom{l}}$ $\gets$ \Call{TransferIntoNormalform}{$t_r$}
		\State \Return $t$
		\EndIf
		\EndFunction
	\end{algorithmic}
\end{algorithm}

\begin{lema} \label{Lemma:equivSubtrees}
	Let $s^n_1, s^n_2, s^n_3, s^n_4$ be dexterity trees and let
	\begin{displaymath}
		{t^{n+1}_1} = \vcenter{\hbox{\begin{tikzpicture}
					\tikzstyle{level 1}=[level distance=7mm,sibling distance=8mm]
					\node {$r$}
					child {node {$s^n_1$}}
					child {node {$s^n_2$}};
		\end{tikzpicture}}}
		, \quad
		{t^{n+1}_2} = \vcenter{\hbox{\begin{tikzpicture}
					\tikzstyle{level 1}=[level distance=7mm,sibling distance=8mm]
					\node {$r$}
					child {node {$s^n_3$}}
					child {node {$s^n_4$}};
	\end{tikzpicture}}}
	\end{displaymath}
	be two dexterity trees with roots with the same label $r \in \{ L, R \}$.
	Then $t^{n+1}_1 \sim t^{n+1}_2$ if and only if 
	$s^n_1 \sim s^n_3$ and $s^n_2 \sim s^n_4$.
\end{lema}
\begin{proof}
	If $s^n_1 \sim s^n_3$ and $s^n_2 \sim s^n_4$ then we can transform
	the subtrees of $t^{n+1}_1$ separately into the corresponding subtrees of $t^{n+1}_2$ and we get $t^{n+1}_1 \sim t^{n+1}_2$.
	
	For the converse direction note that we have the following commutativity for dexterity trees.
	Applying a rule from \eqref{Eq:subtreeChange} to the root and then changes to the left and right subtrees is the same as first applying the subtree changes to the respective other tree and then a suitable rule from \eqref{Eq:subtreeChange} to the root.
	If $t^{n+1}_1 \sim t^{n+1}_2$ 
	there is a sequence of equivalent trees given by applications of \eqref{Eq:subtreeChange}. Since $t^{n+1}_1$ and $t^{n+1}_2$ have a root labelled the same way and by the above argument, we can write this	
	 sequence of equivalent trees as a sequence of changes to the left and right subtree followed by an even number of applications of \eqref{Eq:subtreeChange} to the root. The latter cancel out and we have  a sequence of equivalent subtrees that shows $s^n_1 \sim s^n_3$ and $s^n_2 \sim s^n_4$.
\end{proof}

\begin{lema} \label{Lemma:TwoNormalFormTreesNonEquiv}
	Let $\tau_1^n, \tau_2^n$ be two dexterity trees in normal form with $\tau_1^n \neq \tau_2^n$. 
	Then $\tau_1^n \nsim \tau_2^n$.
\end{lema}
\begin{proof}
	If one of $\tau_1^n$ or $\tau_2^n$ is a fixed tree then no other dexterity tree is equivalent to that one and we are done.
	If both $\tau_1^n$ and $\tau_2^n$ are not fixed trees both of their roots are labelled by $R$. Therefore either their left or right subtrees (or both) differ from each other. Pass to such a pair of subtrees and repeat this process. Since the trees have finite height this terminates and we find two non equivalent subtrees. By \cref{Lemma:equivSubtrees} all of the previous pairs of trees and in particular $\tau_1^n$ and $\tau_2^n$ are non equivalent.
\end{proof}

\begin{proof}[Proof of \cref{Theo:TreeTheorem}]
	By \cref{Lemma:canFindNormalForm} and \cref{Lemma:TwoNormalFormTreesNonEquiv} we see that there is exactly one tree in normal form in each equivalence class in $\treeEqSet{n}$.
	
	Next, we determine $|\treeEqSet{n}|$.
	For $n=1$ the two equivalence classes are `$L$' and `$R$' and we have $|\treeEqSet{1}| =2$.
		
	For $n+1$ with $n \geq 1$ we have seen above that $|\treeEqSet{n+1}|$ is the number of dexterity trees of height $n+1$ that are in normal form. 	
	The normal form trees of height $n+1$ with a root labelled `$L$' are the fixed trees where each pair of siblings is labelled either (`$L$', `$R$') or (`$R$', `$L$'). The number of such trees is $2^{2^{n-1} -1}$.	
	The normal form trees of height $n+1$ with a root labelled `$R$' are of the form where both the left and right subtree of the root is a normal form tree of height $n$. By recursion the number of each of these tress is $|\treeEqSet{n}|$.
	Overall we get the desired recursion formula
	\begin{displaymath}
		|\treeEqSet{n+ 1}| = |\treeEqSet{n}|^2 + 2^{2^{n-1} -1}.
	\end{displaymath}
\end{proof}

\begin{rema}
	The sequence in \cref{Theo:TreeTheorem} can be found in the OEIS as
	\href{http://oeis.org/A332757}{A332757} \cite{oeis}. 
	It also describes  the number of involutions in the n-fold iterated wreath product of $\mathbb{Z} / 2 \mathbb{Z}$ or the number of involutory automorphisms of $\tree {n+1}$.
	The first five terms are 
	$2, 6, 44, 2064$ and  $4292864$.
\end{rema}

\begin{rema}
	It is possible to make choices other than the normal forms for the representatives of $\treeEqSet{n}$. Recall that a normal form tree consists of the constant `$R$' tree with some subtrees swapped with fixed subtrees. But
	one can also choose a different dexterity tree $t^n$ and set the representatives to be of the form $t^n$ with
	some subtrees swapped with fixed subtrees. Algorithm \ref{Alg:TransformIntoNormalform} only needs a small modification which tracks $t^n$ and when it needs to change the label of the current root, to deliver the new representatives.
\end{rema}

\appendix

\section{Some basic lemmas about adjunctions} \label{sec:appAdjoints}

	In this appendix we collect some proofs of statements in \cref{sec:introDual} about well-known properties of adjoints for the reader's convenience.
	

\begin{proof}[Proof of \cref{Lemma:IsomorphismsOfAdjointsAreAdjunctible}]
	The first zig-zag identity \eqref{eq:zig} is given by
	\begin{center}
		\begin{tikzcd}
			&   & l'r'l' \arrow[rd, "\mu^{-1} \times \nu^{-1} \times \id_{l'}", Rightarrow]      &   &    \\
			& l'rl \arrow[ru, "\id_{l'} \times \nu \times \mu", Rightarrow] \arrow[rd, "\mu^{-1} \times \id_{rl}" description, Rightarrow] & & lrl' \arrow[rd, "c \times \id_{l'}", Rightarrow] &    \\
			l' \arrow[ru, "\id_{l'} \times u", Rightarrow] \arrow[rd, "\mu^{-1}"', Rightarrow] &  & lrl \arrow[ru, "\id_{lr} \times \mu" description, Rightarrow] \arrow[rd, "c \times \id_l" description, Rightarrow] &                                                  & l' \\
			& l \arrow[ru, "\id_l \times u" description, Rightarrow]  & & l \arrow[ru, "\mu"', Rightarrow] &   
		\end{tikzcd}
	\end{center}
	where all three squares commute by the interchange law in a bicategory.
	The bottom is equal to $\id_{l'}$ by \eqref{eq:zig} of $(l \adj r, u, c)$.
	Similarly, \eqref{eq:zag} follows from the interchange law and \eqref{eq:zag} of $(l \adj r, u, c)$.
\end{proof}

\begin{proof}[Proof of \cref{Lemma:DifferentUnitsOfSameAdjunction}]
	We define the two $2$-morphisms
	\begin{align*}
		\varphi: \ \big( \quad &r  \ \xRightarrow{ \  \, u \times \id_r \ \, }  \ r \circ l \circ r \
		\xRightarrow{ \ \id_r \times c' \ } \ r \quad \big),
		\\
		\psi: \ \big( \quad &r \ \xRightarrow{ \  u' \times \id_r \ } \ r \circ l \circ r  \
		\xRightarrow{ \ \, \id_r \times c \ \, } \ r \quad \big)  .     
	\end{align*}
	To see that they are inverse to each other we consider the diagram
	\begin{center}
		\begin{tikzcd}
			&                                                                                                           & r \arrow[rd, "u' \times \id_{r}", Rightarrow]                                                                                          &                                                  &   \\
			& rlr \arrow[ru, "\id_{r} \times c'", Rightarrow] \arrow[rd, "u' \times \id_{rlr}" description, Rightarrow] &                                                                                                                                        & rlr \arrow[rd, "\id_{r} \times c", Rightarrow]   &   \\
			r \arrow[ru, "u \times \id_{r}", Rightarrow] \arrow[rd, "u' \times \id_{r}"', Rightarrow] &                                                                                                           & rlrlr \arrow[ru, "\id_{rlr} \times c'" description, Rightarrow] \arrow[rd, "\id_{r} \times c \times \id_{lr}" description, Rightarrow] &                                                  & r \\
			& rlr \arrow[ru, "\id_{rl} \times u \times \id_{r}" description, Rightarrow]                                &                                                                                                                                        & rlr \arrow[ru, "\id_{r} \times c'"', Rightarrow] &  
		\end{tikzcd}
	\end{center}
	where all three squares commute by the
	interchange law in a bicategory.
	The bottom path equals $\id_r$ by \eqref{eq:zig} of the two adjunctions.
	Similarly, we see that $\varphi \circ \psi = \id_r$, where we need \eqref{eq:zag} of the two adjunctions.
	Using the interchange law as above it also follows that
	\begin{displaymath}
		(\psi \times \id_l) \circ u = u' 
		\qquad \text{and} \qquad
		(\id_l \times c) \circ c = c'.
	\end{displaymath}
\end{proof}

\begin{proof}[Proof of \cref{Lemma:composeAdjunctions}]
	The first zig-zag identity \eqref{eq:zig} is given by the diagram
	\begin{center}
		\begin{tikzcd}
			&                                                                                                                                                                                             & f^L g^L g f f^L g^L \arrow[dd, "(ii)", phantom] \arrow[rd, "\id_{f^L} \times c_g \times \id_{f f^L g^L}", Rightarrow] &                                                                                                             &         \\
			& f^L g^L g g^L \arrow[rdd, "\id_{f^L} \times c_g \times \id_{g^L}" description, Rightarrow] \arrow[dd, "(i)", phantom] \arrow[ru, "\id_{f^L g^L g} \times u_f \times \id_{g^L}", Rightarrow] &                                                                                                                       & f^L f f^L g^L \arrow[rdd, "c_f \times \id_{f^L} \times \id_{g^L}", Rightarrow] \arrow[dd, "(iii)", phantom] &         \\
			&                                                                                                                                                                                             & {}                                                                                                                    &                                                                                                             &         \\
			f^L g^L \arrow[rr, "\id", Rightarrow] \arrow[ruu, "\id_{f^L} \times \id_{g^L} \times u_g", Rightarrow] & {}                                                                                                                                                                                          & f^L g^L \arrow[rr, "\id", Rightarrow] \arrow[ruu, "\id_{f^L} \times u_f \times \id_{g^L}" description, Rightarrow]    & {}                                                                                                          & f^L g^L.
		\end{tikzcd}
	\end{center}
	Here $(i)$ commutes as \eqref{eq:zig} of $g^L \adj g$ composed with $(f^L \circ \, - \, )$.
	The triangle $(iii)$ is \eqref{eq:zig} of $f^L \adj f$ composed with $( \, - \,  \circ g^L )$
	and $(ii)$ commutes by the interchange law in a bicategory. 
	The second zig-zag identity \eqref{eq:zag} is similar.
\end{proof}

\section{Dualizability} \label{sec:appDual}

In this Appendix we review some of the definitions of higher dualizability that can be found in the literature and see that they agree with $n$-dualizability as defined in \cref{Def:higherDual}.

Let $\ca C$ be an $\infN N$category.
A $k$-morphism $f$ \emph{has all adjoints} in $\ca C$ 
if there exists a \emph{tower of adjunctions}
\begin{displaymath}
  \dots \, \adj f^{LLL} \adj f^{LL} \adj f^L \adj f \adj f^R \adj f^{RR} \adj f^{RRR} \adj \, \dots
\end{displaymath}
in $\ca C$.
This is the notion of adjunctibility that one gets from Lurie's definition of fully dualizable objects
in \cite{lurie2009classification}. The corresponding data in $\ca C$ can be defined as follows.

\begin{defi}\label{def:fully_n_dualizable}
    Let $\ca C$ be a symmetric monoidal $\infN N$category and
    $X$ an object in $\ca C$.
    A \emph{set of full $1$-dualizability data} for $X$ is a dual object $X^*$ together
    with unit $u$ and counit $c$ witnessing the duality.    
    For $n\geq 3$ a \emph{set of full $n$-dualizability data} is
    a set of full $(n-1)$-dualizability data together with towers of adjunctions for all the unit and counit
    $n$-morphisms in the set of full $(n-1)$-dualizability data.
    An object $X$ is \emph{fully $n$-dualizable} if there exists a set of full $n$-dualizability data
    for $X$.
\end{defi}

Unwrapping this definition for $n \geq 3$ we have the
following data.
 \begin{itemize}
   \item A dual object $X^*$ of $X$ together with unit $u$ and counit $c$,
   \item Towers of adjunctions for $u$ and $c$
   \begin{align}
     \dots \adj u^{LL} \adj u^L \adj \ & u \adj u^R \adj u^{RR} \adj \dots \ , 
     \label{eq:ExampleKis3_1} \\
     \dots \adj c^{LL} \adj c^L \adj \ & c \adj c^R \adj c^{RR} \adj \dots \ ,
     \label{eq:ExampleKis3_2}
   \end{align}
   \item  Towers of adjunctions for all the units and counits in the 
   adjunctions in \eqref{eq:ExampleKis3_1} and \eqref{eq:ExampleKis3_2} as
   well as the corresponding units and counits.
   \item \ldots
   \item Towers of adjunctions for all the unit and counit $n$-morphisms.
 \end{itemize}	    

\begin{rema}
	The above definition is the data in $\ca C$ that one gets from Lurie's definition of
	\emph{fully dualizable objects}. To see this, first discard the non-invertible $k$-morphisms for $k>n$. Therein take  the \emph{fully dualizable} part $\ca C^{\fd}$, which can
	be obtained by discarding $k$-morphisms that are not adjunctible and objects that are not
	dualizable. An object is fully dualizable in Lurie's sense if it lies in the essential image of the natural symmetric monoidal functor
	$\ca C^{\fd} \to \ca C$.
	Furthermore, the $\infN n$category $\ca C^{\fd}$ is said to \emph{have duals}, meaning that every object is fully $n$-dualizable and every morphism has a left and right adjoint.
	See section 2.3 of \cite{lurie2009classification}.
\end{rema}

If $\ca C$ is merely a monoidal $\infN N$category we need to differentiate between left and right duals
for objects.
We adapt the definition of a set of full $n$-dualizability data of an object $X$ to start with
a tower of left and right duals for $X$ of infinite length 
\begin{displaymath}
	\dots \ \adj X^{LL} \adj X^L \adj X \adj X^R \adj X^{RR} \adj \ \dots \ .
\end{displaymath}

In \cref{def:fully_n_dualizable} this tower is hidden behind the braiding.
\begin{lema} \label{Lemma:DeloopingDualizableObject}
	Let $\ca C$ be a symmetric monoidal category. Let $X$ be a dualizable object in $\ca C$ with dual $X^{*}$, unit $u$ and counit $c$.
	Then we have the adjunctions
	\begin{displaymath}
		(X^{*} \adj X, u ,c) \qquad \text{and} \qquad (X \adj X^{*}, u \circ B_{X,X^{*}} ,c \circ B_{X,X^{*}}),
	\end{displaymath}
	in the delooping bicategory $\mathbf{B}\ca C$, where $B$ is the braiding of $\ca C$.
\end{lema}

In \cite{thesisAraujo} Araújo showed that for full $n$-dualizablity one does not need to construct a set of full $n$-dualizability data, but a partial set of data suffices. For bicategories, this statement was already proven in \cite{lurie2009classification} and as Theorem 3.9 in \cite{pstragowski2014dualizable}.
With our definitions we can formulate the slight generalization to not necessarily braided
monoidal categories as follows.

\begin{defi}
    Let $\ca C$ be a monoidal $\infN N$category and
    $X$ an object.
    A \emph{set of partial $1$-dualizability data} for $X$ is a left or right dual object $X^*$ together
    with unit $u$ and counit $c$ witnessing the duality.    
    For $n\geq 2$ a \emph{set of partial $n$-dualizability data} is
    a set of partial $1$-dualizability data together with sets of $a^n$-adjunctibility data for $u$ and
    $c$ for a dexterity map $a^n: \{ 1, 2, \ldots, n \} \to \{ L, R \}$.
\end{defi}

\begin{prop}[Araújo, Theorem 4.1.19 in \cite{thesisAraujo}] \label{Prop:PartialEquivalentToFullDual}
  Let $\ca C$ be a symmetric monoidal $\infN N$category and $n\geq 1$.
  An object $X$ in $\ca C$ is fully $n$-dualizable (\cref{def:fully_n_dualizable}) if and only if it is $n$-dualizable (\cref{Def:higherDual})
  if and only if
  it has
  a set of partial $n$-dualizability data.
\end{prop}
\begin{proof}
	We construct full $n$-dualizability data from a set of partial $n$-dualizability data.
	
	Consider the delooping $\mathbf{B}\ca C$ which is an $(\infty, N+1)$-category with
	a single object and
	$1$-morphisms given by objects in $\ca C$ with composition given by the symmetric monoidal product.
	The partial $n$-dualizability data of $X$ in $\ca C$ gives us $a^{n+1}$-adjunctibility data for $X$ in $\mathbf{B}\ca C$ with a dexterity function
	$a^{n+1}: \{ 1, 2, \ldots, n+1 \} \to \{ L, R \}$.
	By \cref{Lemma:DeloopingDualizableObject} we also have $b^{n+1}$-adjunctibility data for $X$, where $b^{n+1}$ agrees with $a^{n+1}$ except for
	$b^{n+1}(1) = - a^{n+1}(1)$ with the opposite value.
	By \cref{Theo:fromMixedToNadjunctibility} we see that $X$ in $\mathbf{B}\ca C$ is ambidextrous $n$-adjunctible and we can extend the 
	ambidextrous $n$-adjunctibility data to $(n+1)$-adjunctibility data.
	
	By ambidexterity we have towers of adjunctions up to the level of $n$-morphisms.	
	Consider the unit $u$ and counit $c$ of an adjunction of $n$-morphisms $f \adj g$.
	From the $n$-adjunctibility data we have two-step towers
	\[
		u^L \adj u \adj u^R \qquad \text{and} \qquad c^L \adj c \adj c^R .
	\]	
	Applying the Interchange \cref{Lemma:MainLemma} to $f \adj g$ gives us an adjunction $(g \adj f, c^L, u^L)$.
	But we also have an adjunction $g \adj f$ from ambidexterity which has a unit and counit that are part of the $n+1$-adjunctibility data and therefore have left and right adjoints.
	It follows now from \cref{Prop:uniqueness_of_adj} that $u^L$ and $c^L$ have left adjoints $u^{LL}$ and $c^{LL}$.
	Continuing with $(g \adj f, c^L, u^L)$ and two-step towers $c^{LL} \adj c^{L} \adj c$ and $u^{LL} \adj u^{L} \adj u$ we can inductively craft an infinite tower of left adjoints for $u$ and $c$. Analogously, we build towers of right adjoints for $u$ and $c$.
	Altogether we get full $n$-dualizability data for $X$ in $\ca C$.
\end{proof}

\section{Proofs of dualizability in the Morita bicategory} \label{sec:appMorita}

\begin{proof}[Proof of \cref{Prop:allAlgDual}]
  Let $A$ be a $k$-algebra. The identity morphism of $A$ in $\Alg 1$ is given by $_A A_A$.
  We need to show that the zig-zag identities hold up to equivalence.
  For the first one
  \begin{displaymath}
   \left( \  A \ \ \xrightarrow{_k A_{A \te A^{\op}} \, {\te_k} \, _A A_A } \  \ A \te_k A^{\op} \te_k A \
    \xrightarrow{ _A A_A \, {\te_k} \, _{A^{\op} \te A} A _k} \ \ A 
    \ \right) \quad \cong \quad _A A_A,
  \end{displaymath}
  we need a bimodule isomorphism 
  \begin{equation} \label{eq:bimodZigIso}
    \left( _k A_{A \te_k A^{\op}} \ {\te_k} \ _A A_A  \right) 
    \ \bigotimes_{A \te_k A^{\op} \te_k A} \
    \left( _A A_A \ {\te_k} \ _{A^{\op} \te_k A} A _k \right)
    \quad
    \cong \quad _A A_A \, .
  \end{equation}
  This isomorphism is given by mapping an element $(a \te b) \te (c \te d)$ of the left
  hand side in \eqref{eq:bimodZigIso}
  to the product $bdac$.
  Similarly, we get an isomorphism for the other zig-zag identity
  \begin{align*} 
    \left(  A^{\op} 
    \xrightarrow{ {}^{\vphantom{\op}}_{A^{\op}} A_{A^{\op}}^{\op}
     {\te_k} \, _k A_{A \te A^{\op}} } 
     A^{\op} \te_k A \te_k A^{\op} 
    \xrightarrow{ _{A^{\op} \te A} A _k  
     {\te_k} \, {}^{\vphantom{\op}}_{A^{\op}} A_{A^{\op}}^{\op}} \ A^{\op}
     \right) 
    \cong  {}^{\vphantom{\op}}_{A^{\op}} A_{A^{\op}}^{\op}.
  \end{align*}
\end{proof}

Before we turn to the general question, which bimodules $_A M_B$ have adjoints, consider the
special case that $A = B = k$ is a field. Then $_A M_B$ is just a vector space and
it has a left/right adjoint precisely when it has a left/right dual as an object in the category $\Vect$ of vector spaces and linear maps.
To see this, we use a finite basis to construct the unit morphism. We will see that
we  need a finiteness condition for the general case of bimodules to replace the finite basis.
Similar to the vector space case the left and right adjoint of $_A M_B$ will then be given by $\Hom_A(M, A)$ and
$\Hom_B(M, B)$, respectively.

Let $_A L_B$ be a bimodule. Assume we have a right adjoint $_B R_A$ in $\Alg 1$.
Then there is a unit $2$-morphism
\begin{displaymath}
  u: \ {_A A_A} \ \to \ _A L_B \ {\te_B} \ { _B R_A}
\end{displaymath}
of $(A,A)$-bimodules and a counit $2$-morphism
\begin{displaymath}
  c: \ {_B R_A}  \ {\te_A} \ { _A L_B} \ \to \ {_B B_B}
\end{displaymath}
of $(B,B)$-bimodules satisfying the zig-zag identities.
The unit $u$ is completely determined by 
\begin{equation} \label{eq:moduleAdjUn}
  u(1) = \sum_{i,j} l_i {\te_B} r_j \qquad \text{for some} \quad l_i \in L, \ r_j \in R,
\end{equation}
and $c$ is given by a bilinear map
\begin{equation} \label{eq:moduleAdjCo}
  R \times L \to B, \qquad (r,l) \, \mapsto \, c(r,l).
\end{equation}
The first zig-zag identity is 
\begin{equation} \label{eq:moduleAdjZig}
  l \ \xmapsto{ u \te \id } \
  \sum_{i,j} l_i \te r_j \te l 
  \ \xmapsto{ \id \te c } \
   \sum_{i,j} l_i \te c(r_j,l)
  \ \xmapsto{ \sim } \  \sum_{i,j} l_i c(r_j,l) \ = \ l
\end{equation}
for all $l \in L$ and the second zig-zag identity is 
\begin{equation} \label{eq:moduleAdjZag}
  r \ \xmapsto{ \id \te u } \
  \sum_{i,j} r \te l_i \te r_j 
  \ \xmapsto{ c \te \id } \
   \sum_{i,j} c(r,l_i) \te r_j
  \ \xmapsto{ \sim } \  \sum_{i,j} c(r,l_i)r_j \ = \ r
\end{equation}
for all $r \in R$.
Recall that for finite dimensional
vector spaces one uses the fact that there is a finite dual basis to
determine the $l_i$ and $r_j$.
For modules we want finitely presented projectiveness to replace
the finite dimensionality. We will use the fact that a
right $A$-module $P$ is finitely presented projective if and only if 
there is a retract of a free $A$-module $A^n$ of finite rank $n$.
The following Lemma and Proposition and their proofs are mainly taken
from \cite{Yuan}.

\begin{lema}[Dual basis Lemma] \label{Lemma:dualBasis}
  A left $A$-module $P$ is finitely presented projective if and only if 
  there
  exist $e_1, \ldots, e_n \in P$ and elements $e^*_1, \ldots, e^*_n$ 
  in the right $A$-module $P^* := \Hom_A(P,A)$
  such that
  \begin{equation} \label{eq:BasisEq}
    p = \sum_{i = 1}^{n} e_i (e^*_i (p)) \qquad \text{for all} \quad p \in P.
  \end{equation}
  Furthermore, in this case we have
  \begin{equation}   \label{eq:DualBasisEq}
    p^* = \sum_{i = 1}^{n} p^* (e_i) e^*_i \qquad \text{for all} \quad p^* \in P^*.
  \end{equation}   
\end{lema}
\begin{proof}
  The elements $e_1, \ldots, e_n \in P$ define a morphism $r: A^n \to P$
  and the $e^*_1, \ldots, e^*_n : P \to A$ define a morphism
  $\iota: P \to A^n$.
  Equation \eqref{eq:DualBasisEq} implies that $r$ is a retract $r \circ \iota = \id_{P}$.
  Conversely, 
  for a retract $r: A^n \to P$ and map
  $\iota: P \to A^n$ such that $r \circ \iota = \id_{P}$
  we can define $e_1, \ldots, e_n \in P$ as images of the unit vectors in $A^n$
  and $e^*_1, \ldots, e^*_n$ as components of $\iota$.
  
  Equation \eqref{eq:DualBasisEq} follows from applying $\Hom_A( \, \_ \, ,P)$ to get dual morphisms
  $r^*: P^* \to A^n$ and $\iota^*: A^n \to P^*$ with $\iota^* \circ r^* = \id_{P^*}$.
  So $P^*$ is finitely presented projective and the $e^*_i$
  together with the images of $e_i$ under the canonical map $(P^*)^* \to P$ form a
  dual basis.
\end{proof}

\begin{proof}[Proof of \cref{Prop:bimodDual}]
   Let $_A R_B$ be a bimodule that is finitely
   presented projective as a left $A$-module. Define 
   $L = \Hom_A(R, A)$. Then $L$ is $(B,A)$-bimodule.
   By \cref{Lemma:dualBasis}
   there are $e_1, \ldots, e_n \in R$ and $e^*_1, \ldots, e^*_n \in L$ such that
   \begin{equation} \label{eq:DualBasisProp1}
    r = \sum_{i = 1}^{n} e_i (e^*_i (r)) \qquad \text{for all} \quad r \in R.
   \end{equation}
   and
   \begin{equation} \label{eq:DualBasisProp2}
    l = \sum_{i = 1}^{n} l (e_i) e^*_i \qquad \text{for all} \quad l \in L.
    \end{equation}
    We define a unit as in \eqref{eq:moduleAdjUn} with $r_i = e_i$ and
    $l_i = e^*_i$ and a counit as in 
    \eqref{eq:moduleAdjCo} with $c(r,l) = l(r)$ for $r \in R$ and $l \in L$.
   Both zig-zag identities \eqref{eq:moduleAdjZig} and \eqref{eq:moduleAdjZag}
   are satisfied by \eqref{eq:DualBasisProp1} and \eqref{eq:DualBasisProp2}.

   Conversely, let $_A R_B$ be a bimodule with left adjoint $_B L_A$.
   Let $u$ be the unit and $c$ be the counit given by 
   \eqref{eq:moduleAdjUn} and \eqref{eq:moduleAdjCo}.
   The first zig-zag identity is
   equation \eqref{eq:moduleAdjZig} and gives us $ l = \sum_{i,j} l_i c(r_j,l)$ for
   all $l \in L$.
   So $l_1, \ldots, l_n$ together with $r_1, \ldots, r_n$ are a dual basis
   as in \cref{Lemma:dualBasis}. Therefore $_A R_B$ is finitely
   presented projective as a left $A$-module.
   Furthermore $L \cong \Hom_A(R, A)$ by the uniqueness of adjoints up to isomorphism
   from \cref{Lemma:adjointsUnique}.
   
   Dually, the statement for right adjoints can be proven
   analogously or by using the opposite category.
\end{proof}

\section*{Data availability statement}
All data generated or analyzed during this study are contained in this document.

\section*{Conflict of interest statement}
On behalf of all authors, the corresponding author states that there is no conflict of interest.


\bibliographystyle{halpha}
\bibliography{HigherAdjoints_main}

\end{document}